\documentclass[11pt]{article}
\usepackage{}
\usepackage[total={6in, 9in}]{geometry}
\usepackage{amssymb}
\usepackage{mathrsfs}
\usepackage[tbtags]{amsmath}
\usepackage{amsfonts,amsthm,amssymb}
\usepackage{algorithm,algorithmic}
\usepackage{amsfonts}
\usepackage{caption}
\usepackage{float}
\usepackage{graphicx}
\usepackage{graphics}
\usepackage{subfigure}
\usepackage{enumerate}
\usepackage[numbers,sort&compress]{natbib}
\usepackage[
pdfauthor={},
pdftitle={},
pdfstartview=XYZ,
bookmarks=true,
colorlinks=true,
linkcolor=blue,
urlcolor=blue,
citecolor=blue,
bookmarks=true,
linktocpage=true,
hyperindex=true
]{hyperref}
\DeclareMathOperator{\Int}{Int} 
\DeclareMathOperator{\Ints}{int} 
\newtheorem{lem}{Lemma}
\newtheorem{thm}[lem]{Theorem}

\newtheorem{claim}{Claim}

\newtheorem{fact}{Fact}

\newtheorem{Que}[lem]{Question}

\DeclareMathOperator{\df}{def}
\DeclareMathOperator{\odd}{odd}
\newcommand{\CC}{\mathcal{C}}

\def\dist{{\fam0 dist}}

\usepackage[natural]{xcolor}
\usepackage{tikz}
\usepackage{subfigure,tikz}
\usepackage{tikz-3dplot}
\usepackage{pgf,tikz,pgfplots}
\usepackage{mathrsfs}
\usetikzlibrary{shapes.geometric}
\usetikzlibrary{shapes, calc,positioning} 
\begin{document}

\title{2-factors in $\frac{3}{2}$-tough maximal planar graphs}

\author{Lili Hao$^{a}$,\quad Hui Ma$^{a}$,\quad Songling Shan$^{b}$,\quad Weihua Yang$^{a,}$\footnote{Corresponding author. E-mail: ywh222@163.com; yangweihua@tyut.edu.cn}\\
\\ \small $^{a}$Department of Mathematics, Taiyuan University of Technology,\\
\small Taiyuan, Shanxi 030024, China\\
\small $^{b}$Department of Mathematics and Statistics, Auburn University,\\
\small Auburn,  Alabama 36849, U.S.A.\\
}

\date{}
\maketitle

\emph{\textbf{Abstract.}} 
The toughness of a  graph $G$ is defined as the minimum value of $|S|/c(G-S)$ over all cutsets $S$ of $G$ if $G$ is noncomplete, and is defined to be $\infty$ if $G$ is complete.  For a real number $t$, 
we say that $G$ is $t$-tough if its toughness is at least $t$.  
Followed from the classic 1956 result of Tutte, every more than $\frac{3}{2}$-tough planar graph on at least three vertices has a 2-factor. 
In 1999, Owens constructed a sequence of maximal planar graphs with toughness  $\frac{3}{2}-\varepsilon$ for any $\varepsilon >0$, but the graphs do not contain any 2-factor. 
He then posed the question of whether there exists a maximal planar graph with toughness  exactly  $\frac{3}{2}$ and with no 2-factor.  
This question was recently answered affirmatively by the third author. 
This naturally leads to the question: under what conditions does a  $\frac{3}{2}$-tough maximal planar graph contain a 2-factor?
In this paper, we provide a sufficient condition for the existence of 2-factors in $\frac{3}{2}$-tough maximal planar graphs, stated  as a bound on the distance between vertices of degree 3.

\vskip 0.5cm  \emph{\textbf{Keywords.}} 2-factor; plane triangulation; toughness

\section{Introduction}
All graphs discussed are assumed to be simple unless explicitly specified.
Let $G$ be a graph. We denote by $V(G)$ and $E(G)$ respectively the vertex set and edge set of $G$, and we let $n(G)=|V(G)|$, $e(G)=|E(G)|$ and $f(G)$ be the number of faces of $G$ if $G$ is embedded on a surface.  
For $v\in V(G)$, $N_G(v)$ is the set of neighbors of $v$ in $G$. 
For $S\subseteq V(G)$, let $N_G(S)=\left(\bigcup_{x\in S}N_G(x) \right)\setminus S$, and let $G[S]$ and $G-S$ be respectively the subgraph of $G$ induced on $S$ and $V(G)\setminus S$. 
For simplicity,  if $S=\{v\}$ is a singleton, $G-\{v\}$ is written as $G-v$.  
For disjoint subsets $V_1,V_2$ of $V(G)$, we use $E_G(V_1,V_2)$ to denote the set of edges in $G$ with one endvertex in $V_1$ and the other endvertex  in $V_2$, and we let $e_G(V_1,V_2)=|E_G(V_1,V_2)|$. 
If $V_1=\{v\}$ is a singleton, we write $E_G(v,V_2)$ and $e_G(v,V_2)$  respectively for  $E_G(\{v\},V_2)$ and $e_G(\{v\},V_2)$. 
For a subgraph $H$ of $G$ with $V(H)=V_1$, we simply write $E_G(H,V_2)$ and $e_G(H,V_2)$  respectively for  $E_G(V(H),V_2)$ and $e_G(V(H),V_2)$.
For $u,v\in V(G)$, the \emph{distance} between $u$ and $v$ in $G$,  denoted by $\dist_G(u,v)$, is the length of a shortest path connecting $u$ and $v$. 
We say a vertex $v$ of $G$ is \emph{adjacent  to} a subgraph  $H$ of $G$, if $v$ is adjacent in $G$ to a vertex from $V(H)$. 
For any $Z\subseteq V(G)$ and $v\in Z$, we call $v$ is a \emph{$Z$-vertex}.
For two integers $p$ and $q$,  let $[p,q]=\{i\in \mathbb{Z}: p\le i \le q\}$.

Let $c(G)$ denote the number of components of a graph $G$.
The toughness of a graph $G$, denoted by $\tau(G)$, is defined as the minimum value of $|S|/c(G-S)$ over all subsets $S$ of $V(G)$ with $c(G-S) \ge 2$ if $G$ is noncomplete, and is defined to be $\infty$ if $G$ is complete. 
For a real number $t$, 
we say that  $G$ is  $t$-tough if its toughness is at least $t$. 
In 1973, Chv\'{a}tal~\cite{Chvatal1973} introduced the concept of toughness. 
It is easy to see that every cycle is 1-tough and so every graph with a Hamiltonian cycle is 1-tough. 
Conversely, Chv\'{a}tal conjectured that there exists a constant $t_0$ such that every $t_0$-tough graph on at least three vertices is Hamiltonian (Chv\'{a}tal's Toughness Conjecture).
Bauer, Broersma and Veldman~\cite{Bauer2000} have constructed $(\frac{9}{4}-\varepsilon)$-tough graphs that are not Hamiltonian for any $\varepsilon>0$, so $t_0$ must be at least $\frac{9}{4}$ if Chv\'{a}tal's Toughness Conjecture is true.
The conjecture has been verified  for certain classes of graphs including planar graphs, claw-free graphs, co-comparability graphs and chordal graphs. 
For a more comprehensive list of graph classes for which the conjecture holds, see  the survey article by Bauer, Broersma and Schmeichel~\cite{Bauer2006} in 2006.

Chv\'{a}tal's conjecture is known to hold for planar graphs due to Tutte's celebrated 1956 result~\cite{Tutte1956}: every 4-connected planar graph is Hamiltonian.  
This connection arises from the fact that $t$-tough noncomplete graphs are  $\lceil 2t\rceil$-connected.  
Consequently, any noncomplete planar graph with toughness greater than  $\frac{3}{2}$ is 4-connected and so is Hamiltonian. 
Since both $C_3$ and $K_4$ are Hamiltonian, we can then state that every planar graph with toughness more than $\frac{3}{2}$ and at least three vertices is Hamiltonian.  
This leads to the following  natural question:
Is every $\frac{3}{2}$-tough planar graph on at least three vertices Hamiltonian? 

A \emph{2-factor} in a graph $G$ is a spanning 2-regular subgraph. 
Thus a Hamiltonian cycle of $G$ is a 2-factor with only one component.  
In 1999, Owens~\cite{Owens1999} constructed a sequence of maximal planar graphs with toughness  $\frac{3}{2}-\varepsilon$ for any $\varepsilon>0$, but the graphs do not contain even a 2-factor. 
He then posed the question of whether there exists a maximal planar graph with toughness  exactly  $\frac{3}{2}$ and with no 2-factor.  
Bauer, Broersma and Schmeichel in the survey~\cite{Bauer2006} commented that 
``One of the challenging open problems in this area is to determine whether every $\frac{3}{2}$-tough maximal planar graph has a 2-factor. If so, are they all Hamiltonian? We also do not know if a $\frac{3}{2}$-tough planar graph has a 2-factor.'' 
Answering these questions, recently, the third author~\cite{Shan2024} constructed a $\frac{3}{2}$-tough maximal planar graph with no 2-factor. 
In the construction, there are  many pairs of vertices with degree 3 that share a common neighbor. 
In this paper, we aim to find 2-factors in $\frac{3}{2}$-tough maximal planar graphs by restricting the distance between vertices of degree 3. In particular, we obtain the following result.

\begin{thm}\label{thm}
	Let $G$ be a $\frac{3}{2}$-tough maximal planar graph on at least three vertices. 
	Then $G$ has a 2-factor if the distance between any two distinct vertices of degree 3 is at least three. 
\end{thm}

It is interesting to ask whether any $\frac{3}{2}$-tough maximal planar graph on at least three vertices that satisfies the condition as described in Theorem~\ref{thm} is Hamiltonian. 

\begin{Que}
	Let $G$ be a $\frac{3}{2}$-tough maximal planar graph on at least three vertices.  
    Suppose that the distance between any two distinct vertices of degree 3 is at least three.  
	Is $G$ Hamiltonian? 
\end{Que}

\section{Proof of Theorem~\ref{thm}}

Let $S$ and $T$ be two disjoint vertex sets of a graph $G$ and $D$ be a component of $G-(S\cup T)$.
Then $D$ is called an \emph{odd component} (resp. \emph{even component}) if $e_G(D,T)\equiv 1  \pmod{2}$ (resp. $e_G(D,T)\equiv 0 \pmod{2}$). 
For any nonnegative integer $k$, we use $\mathcal{C}_{2k+1}$ to denote the set of all odd components $D$ of  $G-(S\cup T)$ with $e_G(D,T)=2k+1$.
Define $\mathcal{C}=\bigcup_{k\ge 0}\mathcal{C}_{2k+1}$  and $c(S,T)=|\mathcal{C}|$.
Let $\delta(S,T)=2|S|+\sum_{y\in T}d_{G-S}(y)-2|T|-c(S,T)$. 
Note that $\delta(S,T)\equiv 0 \pmod{2}$. We will need Tutte's 2-Factor Theorem, which we state below.

\begin{thm}[Tutte~\cite{Tutte1952}]\label{thmTutte1952}
	A graph $G$ has a 2-factor if and only if $\delta(S,T)\ge 0$ for all $S,T\subseteq V(G)$ with $S\cap T=\emptyset$.
\end{thm}

An ordered pair $(S,T)$ consisting of two disjoint vertex subsets $S$ and $T$ of $G$ is said to be a \emph{barrier} if $\delta(S,T)<0$.
Since $\delta(S,T)\equiv 0 \pmod{2}$, we have $\delta(S,T)\le -2$ if $\delta(S,T)<0$. 
A barrier $(S,T)$ is called a \emph{biased barrier} if $|T|$ is minimum and subject to that $|S|$ is maximum among all the barriers of $G$. 
We will use the following properties of a biased barrier. (These properties were proved by Kanno and the third author in~\cite{Kanno2019} for a barrier $(S,T)$ defined as $|S|$ is maximum and subject to that $|T|$ is minimum among all the barriers of $G$, but the proof works for both definitions.)

\begin{lem}[Kanno and Shan~{\cite[Lemma 3.2]{Kanno2019}}]\label{lemKanno2019}
	Let $G$ be a graph with no 2-factor, and let $(S,T)$ be a biased barrier of $G$. Then the following statements hold:
	\begin{enumerate}[(1)]
		\item $T$ is independent in $G$;
		\item if $D$ is an even component with respect to $(S,T)$, then $e_G(T,D)=0$; 
		\item if $D$ is an odd component with respect to $(S,T)$, then $e_G(y,D)\le 1$ for every $y\in T$;
		\item if $D$ is an odd component with respect to $(S,T)$, then $e_G(x,T)\le 1$ for every $x\in V(D)$. 
	\end{enumerate}
\end{lem}

We will use the operation of edge contractions  in the proof of Theorem~\ref{thm}, which identifies the two endvertices of an edge and removes any resulting loops. 
Contracting a connected subgraph is defined as contracting all the edges in the subgraph iteratively.

Let $G$ be a graph. For $X \subseteq V(G)$, denote by $\mathcal{Y}$ the set of components of $G-X$. 
We define an auxiliary bipartite multigraph $B(X)$ by contracting each component of $G-X$ to a single vertex and deleting edges within $X$. 
In other words, the multigraph $B(X)$ has an edge $xy$ for $x \in X$ and $y\in \mathcal{Y}$ if and only if $x$ has a neighbor in $G$ in the component of $G-X$ corresponding to $y$. 

A \emph{matching} $M$ in a graph $G$ is a set of vertex-disjoint edges.  
A vertex $v$ of $G$ is \emph{covered} or \emph{saturated} by $M$ if $v$ in incident in $G$ with an edge from $M$. 
We denote by $V(M)$ the set of vertices of $G$ that are covered by $M$. 
The matching $M$ is \emph{perfect} if $V(M)=V(G)$. 
For a subset $X \subseteq V(G)$, we denote by $\odd(G-X)$ the number of components of $G-X$ that are of odd order. 
Define $\df(G)=\max\{\odd(G-X)-|X|: X\subseteq V(G)\}$. 
Applying Tutte's Theorem~\cite{Tutte1947} to the graph obtained from $G$ by adding $\df(G)$ vertices and joining all edges between $V(G)$ and these new vertices, Berge~\cite{Berge1958} observed that the maximum size of a matching in $G$ is $\frac{1}{2}(n(G)-\df(G))$. 
West~\cite{West2011} gave a short proof of this result and in particular, he proved the following result.

\begin{lem}[West~{\cite[Lemmas 2 and 3]{West2011}}]\label{lemWest2011}
	Let $G$ be a multigraph and $X\subseteq V(G)$ be maximal with $\odd(G-X)-|X|=\df(G)$.  Then the following statements hold:
	\begin{enumerate}[(1)]
		\item every component of $G-X$ has odd order;
		\item for every component $D$ of $G-X$ and any $u \in V(D)$, $D-u$ has a perfect matching; 
		\item the bipartite multigraph $B(X)$ has a matching covering $X$.
	\end{enumerate}
\end{lem}

\proof[\textbf{Proof of Theorem~\ref{thm}}]
Let $G$ be a $\frac{3}{2}$-tough maximal planar graph on at least three vertices such that the distance between any two distinct vertices of degree 3 is at least three. 
Suppose to the contrary that $G$ does not have any 2-factor. 
Thus $G$ is not a complete graph and so is 3-connected by $\tau(G)\ge \frac{3}{2}$. 
This further implies that $|V(G)| \ge 7$ ($G$ would otherwise have a Hamiltonian cycle by Dirac's Theorem). 
We embed $G$ in the plane, and still use $G$ to denote the embedding.   
As $G$ is maximal planar, the embedding is a plane triangulation. 

By Theorem~\ref{thmTutte1952}, $G$ has a barrier.
We let $(S,T)$ be a biased barrier of $G$. Then $(S,T)$ satisfies all the properties as listed in Lemma~\ref{lemKanno2019}.
Let  $ U=V(G)\setminus (S\cup T)$, and for each integer $k\ge 0$, let 
\begin{eqnarray*}
    c_{2k+1} &=& |\mathcal{C}_{2k+1}|, \quad \text{and} \\ 
	U_{2k+1} &=& \bigcup_{D\in \CC_{2k+1}} V(D). 
\end{eqnarray*}

For each $y\in T$, let
\begin{eqnarray*}
    c^{\prime}(y) & = & |\{D\in \mathcal{C}: e_G(y,D)=1\}|,\quad \text{and}\\
    c(y) & = & |\{D\in \mathcal{C}\setminus \mathcal{C}_1: e_G(y,D)=1\}|.
\end{eqnarray*}
By Lemma~\ref{lemKanno2019}(1)-(3), we have $e_G(y,U)=c^{\prime}(y)$ and $d_G(y)=e_G(y,S)+c^{\prime}(y)$.

As $\delta(S,T)\le -2$ and $\sum_{v\in T}d_{G-S}(v)\ge c(S,T)$, it follows that $T\neq \emptyset$. 
Since $|T|$ is minimum among all the barriers of $G$, $(S,T\setminus \{y\})$ is not a barrier for any $y\in T$. 
Thus we have
\begin{eqnarray*}
0\le \delta (S,T\setminus \{y\})
    & = & 2|S|-2|T|+2+\sum_{v\in T\setminus \{y\}}d_{G-S}(v)-c(S,T\setminus \{y\})\\
    & \le & 2|S|-2|T|+2+\sum_{v\in T}d_{G-S}(v)-c^{\prime}(y)-\left (c(S,T)-c^{\prime}(y)\right)\\
    & = & \delta(S,T)+2\le 0.
\end{eqnarray*}
Hence $\delta(S,T)+2=0$. 
On the other hand, by Lemma~\ref{lemKanno2019}, we have 
\begin{equation*}
    \sum_{v\in T}d_{G-S}(v)=\sum_{k\ge 0}(2k+1)c_{2k+1}.
\end{equation*}
As $c(S,T)=\sum_{k\ge 0}c_{2k+1}$, 
we then have 
\begin{equation*}
    \delta(S,T)+2=2|S|-2|T|+\sum_{k\ge 1}2kc_{2k+1}+2=0.
\end{equation*}
This gives
\begin{equation}\label{eq:T-size}
    |T|=|S|+\sum_{k\ge 1}kc_{2k+1}+1.
\end{equation}

Define
\begin{equation*}
\begin{split}
    T_0 & = \{y\in T: c(y)=0\},\quad T_1=\{y\in T: c(y)=1\},\\
    T_2 & = T\setminus (T_0\cup T_1),\ \ \ \ \ \ \ \ \ \quad p=\sum_{v\in T_2}(c(v)-2).
\end{split}
\end{equation*}
Note that 
\begin{equation}\label{eq:T2-size}
    |T_2|=\frac{1}{2}\left(\sum_{k\ge 1}(2k+1)c_{2k+1}-|T_1|-p \right).
\end{equation}

Our proof proceeds in two main stages. 
In the first stage, we establish an upper bound on $c_3$, since the triangles in $\mathcal{C}_3$ pose challenges for constructing an effective cutset. 
In the second stage, we derive a lower bound on the size of $S$. 
Using these bounds, we construct a cutset $W$ in $G$ and show that $|W|/c(G-W)<\frac{3}{2}$ by applying the two inequalities.

\smallskip 

{\bf \noindent Stage 1: Establish an upper bound on $c_3$.}

\smallskip 
We construct an auxiliary graph $H$ by applying the following operations on $G$. 
Roughly, the vertices of $H$ correspond to components of $\mathcal{C}\setminus \mathcal{C}_1$, and an edge $uv$ of $H$ corresponds to a path of length two in $G$ that connects the two components of  $\mathcal{C}\setminus \mathcal{C}_1$ that respectively correspond to  $u$ and $v$. 
We will use the structure of $H$ in deriving  an upper bound on $c_3$.

\noindent{\bf Operations:}
\begin{enumerate}[(1)]
    \item We delete all odd components in $\mathcal{C}_1$, and we map each odd component $D$ of $\mathcal{C}\setminus \mathcal{C}_1$ into a vertex $h(D)$ by contracting $D$. 
        (We will use $h^{-1}(h(D))$ to denote the component $D$.)
        Call the resulting multigraph $G_1$, which is still plane.

    \item In $G_1$, delete all vertices in $S$ and in $T_0$. 
        Call the resulting multigraph $G_2$, which is also plane. 
        Let $\mathcal{C}_{\ge 5}=\bigcup_{k\ge 2}\mathcal{C}_{2k+1}$,  
        \begin{equation*}
        h(\mathcal{C}_3)=\{x\in V(G_2)\setminus T:h^{-1}(x)\in \mathcal{C}_3\},\quad \text{and}\quad 
        h(\mathcal{C}_{\ge 5})=\{x\in V(G_2)\setminus T:h^{-1}(x)\in \mathcal{C}_{\ge5}\}.
        \end{equation*}

    \item In $G_2$, we ``smooth'' each vertex of degree 2. 
        That is, if $v$ is of degree 2 in $G_2$ with neighbors $u$ and $w$, we replace the path $uvw$ with the edge $uw$. 
        As each vertex $x\in h(\mathcal{C}_3) \cup h(\mathcal{C}_{\ge 5})$ satisfies $d_{G_2}(x)\ge 3$, those smoothed vertices are from $V(G_2)\cap T_2$. 
        Then we delete all vertices from $T_1$. 
        Call the resulting multigraph $G_3$, which is still plane. 
        Note that $V(G_3)$ consists of two types of vertices: vertices contracted from odd components of $G-(S\cup T)$ and vertices $y\in T_2$ with $c(y)\ge 3$.

    \item In $G_3$, we split each vertex of $T_2\cup h(\mathcal{C}_{\ge 5})$ which has degree at least 4 in $G_3$ into some independent vertices as follows.
    
    \begin{itemize}
    	\item 
        For each vertex $y\in T_2$ with $d_{G_3}(y)\ge 4$, we split $y$ into $\left\lfloor     \frac{d_{G_3}(y)-2}{2}\right\rfloor$ independent vertices each with degree 2, and one with degree $d_{G_3}(y)-2\left\lfloor \frac{d_{G_3}(y)-2}{2}\right\rfloor$.  
        The last vertex is also of degree 2 if $d_{G_3}(y)$ is even and is of degree 3 otherwise.
    	\item 
        Let $\ell =\sum_{k\ge 2}(2k+1)c_{2k+1}$. 
        Then we have $\sum_{x\in h(\mathcal{C}_{\ge 5})}d_{G_3}(x)\le \ell$. 
        We identify all vertices of $ h(\mathcal{C}_{\ge 5})$ in $G_3$ into a single vertex and split the identified vertex into 
        $$\left\lceil \frac{1}{3} \left(\sum_{x\in  h(\mathcal{C}_{\ge 5})} d_{G_3}(x) \right)   \right\rceil \le \frac{1}{3} (\ell +2)$$ 
        independent vertices.  
    \end{itemize}
   
        Denote the multigraph resulting from $G_3$ by $H$. 
        Note that $\Delta(H)\le 3$. 
        For a vertex $x\in T_2$ which was split in Case 1 above, we denote by $\mathcal{S}(x)$ the set of all vertices of $H$ that are split from $x$. 
        In Case 2, we denote by $\mathcal{S}(h(\mathcal{C}_{\ge 5}))$ the set of all vertices of $H$ that are split from the single identified vertex of all vertices of $h(\mathcal{C}_{\ge 5})$. 
        In both cases, we call each of these vertices a \emph{splitting vertex}. 
\end{enumerate} 

Note that the multigraph $H$ may no longer  be a plane graph  by our operation of identifying vertices of $h(\mathcal{C}_{\ge 5})$. 
However, if a component of $H$ does not contain any splitting vertex from $\mathcal{S}(h(\mathcal{C}_{\ge 5}))$, then that component is still a plane graph.

For a connected subgraph $D$ of $H$ with $V(D) \subseteq h(\mathcal{C}_3)$, we let $h_1^{-1}(D) $ denote the subgraph of $G$  obtained from $D$ by reversing Operations (1) and (3). 
That is, the graph obtained from $D$ by replacing each vertex $x\in V(D)$ by $h^{-1}(x)$, and by replacing each edge $xy \in E(D)$ first by $xwy$, where $w$ was the degree 2 vertex smoothed out in Operation (3) in getting $xy$, then by replacing $xwy$ with $x^{\ast}wy^{\ast}$ such that $x^{\ast}$ and $y^{\ast}$ are respectively the neighbor of $w$ in $G$ from $h^{-1}(x)$ and $h^{-1}(y)$.  
If $W$ is a facial walk of $h_1^{-1}(D)$, then we denote by $\Int(W)$ the subgraph of $G$ induced by $V(W)$ together with all the vertices of $G$ that are embedded in the interior of $W$, and let   $\Ints(W)=\Int(W)-V(W)$.   
Since each edge of $D$ corresponds to a path of length two in $G$,  and $G$ is a plane triangulation, it follows that $\Ints(W)$ contains at  least one vertex from $S$.



We will establish an upper bound on $c_3$ by constructing a matching of $H$ with  appropriate size.  To this end, we first prove the following claim.

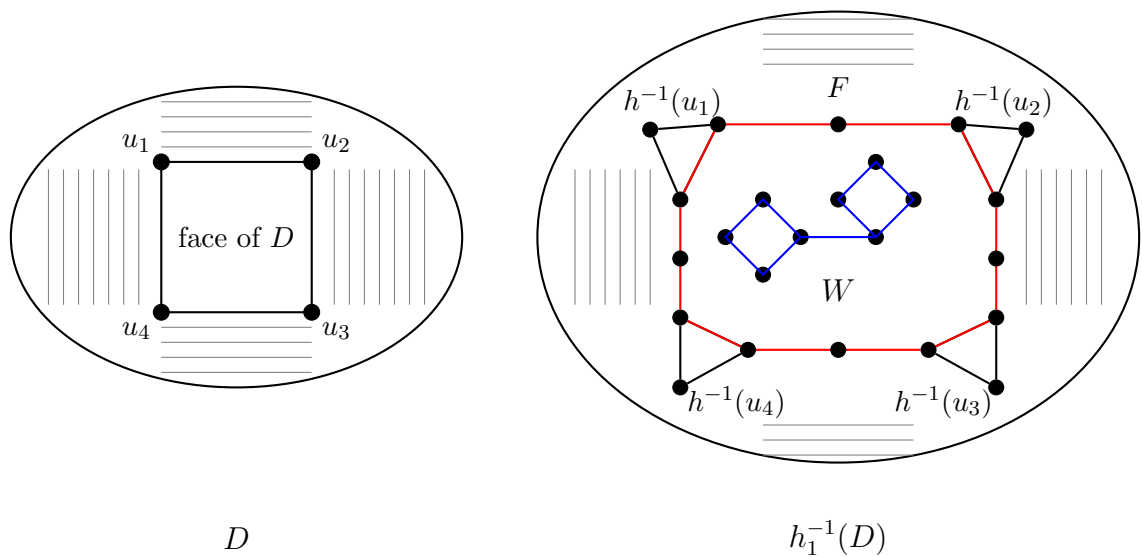
\begin{figure}[!htb]
	\begin{center}
\begin{tikzpicture}
	\draw[thick] (0,0) ellipse (3 and 2);
	\node at (0,-4) {\large $D$};
	
	\draw[thick] (-1,1) -- (1,1) -- (1,-1) -- (-1,-1) -- cycle;
	\node at (0,0) {face of $D$};
	
	\filldraw (-1,1) circle (3pt) node[above left] {$u_1$};
	\filldraw (1,1) circle (3pt) node[above right] {$u_2$};
	\filldraw (1,-1) circle (3pt) node[below right] {$u_3$};
	\filldraw (-1,-1) circle (3pt) node[below left] {$u_4$};
	
	\foreach \x in { -2.5, -2.3, ...,-1.2} {
		\draw[gray] (\x,-0.9) -- (\x,0.9);
	}

	\foreach \x in { 2.5, 2.3, ...,1.2} {
		\draw[gray] (\x,-0.9) -- (\x,0.9);
	}

	\foreach \y in { 1.2,1.4,...,1.8} {
		\draw[gray] (-1,\y) -- (1,\y);
	}

	\foreach \y in { -1.2,-1.4,...,-1.8} {
		\draw[gray] (-1,\y) -- (1,\y);
	}
	
    \begin{scope}[shift={(8,0)}]
	\draw[thick] (0,0) ellipse (4 and 3);
	
    \node[circle, draw, fill=black, inner sep=2pt] (u11) at (-1.6,1.5) {};
    \node[circle, draw, fill=black, inner sep=2pt] (u12) at (-2.1,0.5) {};
    \node[circle, draw, fill=black, inner sep=2pt] (u13) at (-2.5,1.43) {};
		
    \draw[thick] (u11) -- (u12) -- (u13)--(u11);
	\draw[thick,red] (u11) -- (u12); 
 
    \node[circle, draw, fill=black, inner sep=2pt] (u21) at (1.6,1.5) {};
    \node[circle, draw, fill=black, inner sep=2pt] (u22) at (2.1,0.5) {};
    \node[circle, draw, fill=black, inner sep=2pt] (u23) at (2.5,1.43) {};
	 
    \draw[thick] (u21) -- (u22) -- (u23)--(u21);

    \node[circle, draw, fill=black, inner sep=2pt] (u31) at (-1.2,1.5-3) {};
    \node[circle, draw, fill=black, inner sep=2pt] (u32) at (-2.1,0.5-2.5) {};
    \node[circle, draw, fill=black, inner sep=2pt] (u33) at (-2.1,1.43-2.5) {};
	 
    \draw[thick,red] (u21) -- (u22);

    \draw[thick] (u31) -- (u32) -- (u33)--(u31);
    \draw[thick,red] (u31) -- (u33);

    \node[circle, draw, fill=black, inner sep=2pt] (u41) at (1.2,1.5-3) {};
    \node[circle, draw, fill=black, inner sep=2pt] (u42) at (2.1,0.5-2.5) {};
    \node[circle, draw, fill=black, inner sep=2pt] (u43) at (2.1,1.43-2.5) {};

	\draw[thick] (u41) -- (u42) -- (u43)--(u41);
	 
	\draw[thick,red] (u41) -- (u43); 
	
	\draw[thick,red] (u11) -- (u21);
	\path (u11) -- node[midway, fill=black, circle, inner sep=2pt, draw] (M) {} (u21);

	\draw[thick,red] (u12) -- (u33);
    \path (u12) -- node[midway, fill=black, circle, inner sep=2pt, draw] (M) {} (u33);

	\draw[thick,red] (u31) -- (u41);
	\path (u31) -- node[midway, fill=black, circle, inner sep=2pt, draw] (M) {} (u41);
	 
	\draw[thick,red] (u22) -- (u43);
	\path (u22) -- node[midway, fill=black, circle, inner sep=2pt, draw] (M) {} (u43);
	 
	\node at (-2.2,1.8) {$h^{-1}(u_1)$};
	\node at (2.2,1.8) {$h^{-1}(u_2)$};
	\node at (-1.35,-2.2) {$h^{-1}(u_4)$};
	\node at (1.4,-2.2) {$h^{-1}(u_3)$};
	\node at (0,2) {$F$};
	
	
	\node[circle, draw, fill=black, inner sep=2pt] (a1) at (-1.5,0) {};
	\node[circle, draw, fill=black, inner sep=2pt] (a2) at (-1,0.5) {};
	\node[circle, draw, fill=black, inner sep=2pt] (a3) at (-0.5,0) {};
	\node[circle, draw, fill=black, inner sep=2pt] (a4) at (-1,-0.5) {};
	
	\node[circle, draw, fill=black, inner sep=2pt] (b1) at (0,0.5) {};
	\node[circle, draw, fill=black, inner sep=2pt] (b2) at (0.5,1) {};
	\node[circle, draw, fill=black, inner sep=2pt] (b3) at (1,0.5) {};
	\node[circle, draw, fill=black, inner sep=2pt] (b4) at (0.5,0) {};
	
	\draw[thick, blue] (-1.5,0) -- (-1,0.5) -- (-0.5,0) -- (-1,-0.5) -- cycle;
	\draw[thick, blue] (0,0.5) -- (0.5,1) -- (1,0.5) -- (0.5,0) -- cycle;
    \draw[thick, blue] (-0.5,0) --(0.5,0); 
	\node at (0,-0.7) {\large $W$}; 
	 
	\node at (0,-4) {\large $h_1^{-1}(D)$};

	\foreach \x in { -3.5, -3.3, ...,-2.5} {
		\draw[gray] (\x,-0.9) -- (\x,0.9);
	}

	\foreach \x in { 3.5, 3.3, ...,2.5} {
		\draw[gray] (\x,-0.9) -- (\x,0.9);
	}

	\foreach \y in { 2.3,2.5,...,3} {
		\draw[gray] (-1,\y) -- (1,\y);
	}

	\foreach \y in { -2.5,-2.7,...,-3} {
	\draw[gray] (-1,\y) -- (1,\y);
}
	
	\end{scope}
\end{tikzpicture}
\end{center}
	\caption{A depeiction of $D$, $h_1^{-1}(D)$, $F$, and $W$ for Claim~\ref{claim:componentH-X}.}
\label{f0}
\end{figure}

\begin{claim}\label{claim:componentH-X}
    Let $X\subseteq V(H)$ be maximal with $\odd(H-X)-|X|=\df(H)$, and let $D$ be a component of  $H-X$ with the following properties:
    \begin{enumerate}[(1)]
        \item $|V(D)|\ge 3$ is odd;
        \item $V(D)\subseteq h(\mathcal{C}_3)$;
        \item $e(D)\ge 1.5n(D)-0.5$. 
    \end{enumerate}
    Then $h_1^{-1}(D)$ has at least two faces $F$ such that $F$ satisfies the following properties:  
    for the outer facial walk $W$ of $\Ints(F)$ (here we take the face of $\Ints(F)$ that has an edge joining to a vertex of $F$ as the outer face of $\Ints(F)$), we have $V(W)\subseteq S \cup U_1$, and each $S$-vertex of $W$ is contained in a cycle of $W$. 
    See Figure~\ref{f0} for an illustration of $W$.    
\end{claim}

\proof 
Let $D$ be a component of $H-X$ with the three properties stated above. 
Note that $D$ is a subgraph of $G_3$ containing no splitting vertex, and so every edge $uv$ of $D$ corresponds to a vertex $y\in T_2$ such that $c(y)=2$. 
As $D$ contains no splitting vertex of $G_3$, $D$ is a plane graph. 
By Lemma~\ref{lemKanno2019}(3)-(4), $D$ is simple. 
Recall that $n(D)$, $e(D)$, and $f(D)$, respectively, denotes the number of vertices, edges, and faces of $D$. 
By Euler's formula, we have $f(D)=e(D)-n(D)+2$.  
Since $e(D)\ge 1.5n(D)-0.5$, $e(D) \le 1.5 n(G)$, and every edge of $D$ corresponds to a vertex $y\in T_2$ such that $c(y)=2$, it follows that $e_G(V(h_1^{-1}(D)) \cap U,T) \le 1$.  
Thus, there is at most one facial walk $F$ of $h_1^{-1}(D)$ such that $e_G(V(h_1^{-1}(F)),T_1) =1$, which we call it \emph{bad}; a facial walk of $h_1^{-1}(D)$ that is not bad is \emph{good}. 
Then $h_1^{-1}(D)$ has at least $f(D)-1$ good facial walks.

Let $F$ be a good facial walk of $h_1^{-1}(D)$. 
Suppose $\Ints(F)$ contains an $S$-vertex $x$ that satisfies Property $({\ast})$ as listed below:
\begin{enumerate}[(a)]
 	\item $x$ is not adjacent in $G$ to any component of $\mathcal{C} \setminus \mathcal{C}_1$ that is not corresponding to a vertex of $D$,  
 	\item $x$ is not adjacent in $G$ to any component $R$ of $\mathcal{C}_1$ such that $e_G(R, V(h_1^{-1}(D))\cap T) =0$, and 
 	\item $x$ is not adjacent in $G$ to any vertex of $T\setminus V(h_1^{-1}(D))$.  
\end{enumerate}
Then we select such an $S$-vertex $x$. 
Let $S_D$ be the set of all these selected $S$-vertices range over all good facial walks of $h_1^{-1}(D)$.

Consider first that $|S_D|\ge f(D)-2$. 
Each vertex of $D$ is corresponding in $G$ to a component from $\mathcal{C}_3$ of $G-(S\cup T)$. 
Let $G^1_D$ be the union of all components from $\mathcal{C}_3$ that each corresponds to a vertex of $D$, and $T_D=\{y\in T: c(y)=e_G(y,G^1_D)=2\}$. 
Note that each vertex from $T_D$ corresponds to exactly one edge of $D$, and so $|T_D|=|E(D)|$. 
Let $G_D=G[V(G^1_D)\cup T_D]$. 
Note that $G_D=h_1^{-1}(D)$. 
For notation simplicity, let $\mathcal{C}_D$ be the set of components from $\mathcal{C}$ that are contained in $G_D$. 
Then by the construction of $H$, we have
\begin{equation*}
	|T_D|=e(D), \quad \text{and}\quad |\mathcal{C}_D|=n(D).
\end{equation*}
Let $S^{\ast}=S\setminus S_D$ and $T^{\ast}=T\setminus T_D$.  
Since each vertex of $S_D$ satisfies Property $(\ast)$ and each vertex $y\in T_D$ satisfies $c(y)=e_G(y,G^1_D)=2$, 
we have 
\begin{eqnarray*}
	\delta(S^{\ast},T^{\ast})
		& = & 2|S|-2|S_D|-2|T|+2|T_D|+\left(\sum_{y\in T}d_{G-S}(y)-2e(D)-\sum_{y\in T_D,D_0\in \mathcal{C}_1}e_G(y,D_0)\right)\\
		&& \quad -\left(c(S,T)-|\mathcal{C}_D|-\sum_{y\in T_D,D_0\in \mathcal{C}_1}e_G(y,D_0)\right)\\
		& \le & 2|S|-2(f(D)-2)-2|T|+\sum_{y\in T}d_{G-S}(y)-c(S,T)+n(D)\\
		& = & 2|S|-2|T|+\sum_{y\in T}d_{G-S}(y)-c(S,T)-2e(D)+3n(D) \\
		&& \quad \text{(we used $n(D) -e(D)+f(D) =2$ in getting the equation above)}\\
		& \le & 2|S|-2|T|+\sum_{y\in T}d_{G-S}(y)-c(S,T)+1\\
		& = & \delta(S,T)+1<0.
\end{eqnarray*}
Thus $(S^{\ast},T^{\ast})$ is also a barrier of $G$. 
Since $|T_D|=e(D)\ge 1.5n(D)-0.5$ and $n(D)\ge 3$ by the assumptions on $D$, we have $T_D \neq \emptyset$. 
Thus $(S^{\ast},T^{\ast})$ is a barrier of $G$ with $|T^{\ast}|<|T|$, a contradiction to the choice of $(S,T)$.

Thus $|S_D| \le f(D)-3$. This implies that $h_1^{-1}(D)$ has at least two good facial walks $F$ such that every $S$-vertex of the outer facial walk $W$ of $\Ints(F)$ is adjacent in $G$ to some  component of $\mathcal{C}\setminus \mathcal{C}_1$ that is not corresponding to any vertex of $D$, or is adjacent in $G$ to a component $R$ of $\mathcal{C}_1$ such that $e_G(R, V(h_1^{-1}(D))\cap T)=0$, or is adjacent in $G$ to a vertex of $T\setminus V(h_1^{-1}(D))$. 
As $G$ is a plane triangulation and $F$ is a good facial walk of $h_1^{-1}(D)$, it follows that  every $S$-vertex of $W$ is contained in a cycle $C$ of $W$ such that $V(C)$ is a cutset of $G$. 
We have $V(W)\subseteq S \cup U_1$ as $F$ is a good facial walk of $h_1^{-1}(D)$. 
\qed

A component $D$ of $H-X$ that satisfies all properties in Claim~\ref{claim:componentH-X} is called a \emph{compact-component}, and a walk $W$ guaranteed by Claim~\ref{claim:componentH-X} is called an \emph{auxiliary walk} of $D$. 
We let $q$ be the total number of compact-components of $H-X$. 
By Claim~\ref{claim:componentH-X}, every compact-component has at least two auxiliary walks.

\begin{fact}\label{fact:fact1}
    For each component $D$ of $H-X$ that is not a compact-component, by Claim~\ref{claim:componentH-X}, $D$ satisfies one of the following properties:
    \begin{enumerate}[(1)]
        \item $D$ contains a  splitting vertex;
		\item $|V(D)|=1$ and $V(D)\subseteq  h(\mathcal{C}_3)$;
        \item $|V(D)|\ge 3$, $V(D)\subseteq  h(\mathcal{C}_3)$ and $e(D) \le  1.5n(D)-1.5$ (note that $|V(D)|$ is odd by Lemma~\ref{lemWest2011}(1)). 
    \end{enumerate}     
    For $D$ that satisfies (2) or (3) above, we have 
    \begin{equation}
		e_{G_1}(D,T_1)\ge (3n(D)-e_H(D,X))-(3n(D)-3)=3-e_H(D,X),  \label{eq:D-T1-edge}
	\end{equation}
    where recall that $G_1$ was obtained from $G$ by deleting each component in $\mathcal{C}_1$ of $G-(S\cup T)$, and then contracting each of the rest components of $G-(S\cup T)$ into a single vertex.  
\end{fact}

By Lemma~\ref{lemWest2011}(3), the bipartite multigraph $B(X)$ has a matching $M_1$ covering $X$.  
For easier description, we will also  refer a component of $H-X$ as a vertex of $B(X)$. 
For each $i\in [0,2]$, we let
$$
    \mathcal{Y}_i = \{D: \text{$D$ is a noncompact-component of $H-X$ with $V(D) \subseteq h(\mathcal{C}_3)$ and $e_H(D,X)=i$}\},
$$
and let 
$$
    \mathcal{Y}_3 = \{D: \text{$D$ is a noncompact-component of $H-X$ with $V(D) \subseteq h(\mathcal{C}_3)$ and $e_H(D,X)\ge 3$}\}.
$$

We next claim that $B(X)$ has a matching  $M_2$ covering  $\mathcal{Y}_3$ and  at least $\left\lceil \frac{1}{3} (|\mathcal{Y}_1|+2|\mathcal{Y}_2|) \right\rceil$ vertices of $\mathcal{Y}_1 \cup \mathcal{Y}_2$.  
To prove this, we let $B^{\ast}(X)$ be obtained from $B(X)$ by adding a set $Z$ of $\left\lfloor \frac{1}{3} (2|\mathcal{Y}_1|+|\mathcal{Y}_2|)\right \rfloor$ new vertices, and adding edges between $Z$ and $\mathcal{Y}_1 \cup \mathcal{Y}_2$ so that each vertex of $Z$ gets degree 3, and each vertex of $\mathcal{Y}_1 \cup \mathcal{Y}_2$  gets degree at most 3 in the resulting multigraph $B^{\ast}(X)$. 
As $|Z|=\left\lfloor \frac{1}{3}(2|\mathcal{Y}_1|+|\mathcal{Y}_2|)\right \rfloor$, it follows that all but at most one vertex of $\mathcal{Y}_1 \cup \mathcal{Y}_2$ have degree 3 in $B^{\ast}(X)$, and the exceptional vertex, whose existence depends on the remainder of $2|\mathcal{Y}_1|+|\mathcal{Y}_2|$ when divided by three, has degree 1 or 2.   
It suffices to show that $B^{\ast}(X)$ has a matching covering $\mathcal{Y}_1 \cup \mathcal{Y}_2 \cup \mathcal{Y}_3$.  
Suppose not. 
Then there exists $A \subseteq \mathcal{Y}_1 \cup \mathcal{Y}_2 \cup \mathcal{Y}_3$ such that $|N_{B^{\ast}(X)}(A)|<|A|$. 
Let $x\in A$ such that $d_{B^{\ast}(X)}(x)$ is smallest among the degrees of all vertices of $\mathcal{Y}_1 \cup \mathcal{Y}_2 \cup \mathcal{Y}_3$ in $B^{\ast}(X)$.  
Note that $d_{B^{\ast}(X)}(x) \ge 1$ by our comment above. 
Now we get 
\begin{eqnarray*}
    d_{B^{\ast}(X)}(x)+3(|A|-1)\le e_{B^{\ast}(X)}(A,N_{B^{\ast}(X)}(A))\le 3|N_{B^{\ast}(X)}(A)|\le 3(|A|-1), 
\end{eqnarray*}
a contradiction.  
Thus $B^{\ast}(X)$ has a matching $M^{\ast}$ covering  $\mathcal{Y}_1\cup \mathcal{Y}_2\cup \mathcal{Y}_3$. 
Deleting from $M^{\ast}$ edges incident with vertices of $Z$ gives us a desired matching $M_2$.

Let $Y$ be the set of vertices of $\mathcal{Y}_1 \cup \mathcal{Y}_2 \cup \mathcal{Y}_3$ that are covered by $M_2$. 
We claim next that we can construct a matching $M$ of $B(X)$ that covers both $X$ and $Y$. 
Let $F$ be a component of the graph formed by the union of edges from $M_1$ and $M_2$.
Then $F$ can be an even cycle, a path of even length, or a path of odd length.
If $F$ is an even cycle or a path of odd length, then $F$ has a matching $M_F$ covering $V(F)$.
If $F$ is a path of even length and start and so end at two vertices from $X$, we let $M_F=E(F)\cap M_1$ be a matching of $F$ (this case actually does not exist as $M_1$ covers all vertices from $X$).
If $F$ is a path of even length and start and so end at two vertices from $Y$, we let $M_F=E(F)\cap M_2$ be a matching of $F$ (this case similarly does not exist since $M_2$ covers all vertices from $Y$).
If $F$ is a path of even length and start and end at two vertices from $V(B(X)) \setminus (X\cup Y)$ and $Y$, respectively, we let $M_F=E(F)\cap M_2$ be a matching of $F$.
In all the cases above, by the construction, $M_F$ saturates all vertices in $(X\cap V(F))\cup (Y\cap V(F))$. 
Let $M$ be the union of matchings $M_F$ for all the components $F$. 
By the construction, $M$ is a matching with the desired property. 
Thus in $H$, there is a matching $M$ covering $X$, $\mathcal{Y}_3$, and at least $\left\lceil \frac{1}{3} (|\mathcal{Y}_1|+2|\mathcal{Y}_2|)\right \rceil$ vertices of $\mathcal{Y}_1 \cup \mathcal{Y}_2$. 
Note that $|V(M)\cap V(D)| \le 1$ for each component $D$ of $H-X$ and $|V(M)\cap V(D)|=1$ if $D\in \mathcal{Y}_3 $.

For an edge $xy \in M$ with $x\in X$ and $y$ corresponding to a component $D_y$ of $H-X$, we let $y^{\ast}\in V(D_y)$ be a vertex such that $xy^{\ast}\in E(H)$, and we let $M^{\prime}=\{xy^{\ast}: xy \in M\}$.  
Let $D$ be a component of $H-X$.  
If $|V(D)\cap V(M^{\prime})|=1$, we let $y^{\ast}\in V(D)\cap V(M^{\prime})$, and let $M_D$ be a perfect matching of $D-y^{\ast}$ by Lemma~\ref{lemWest2011}(2). 
If $D$ is not covered by $M$ and $D$ contains a splitting vertex $x$, we again let $M_D$ be a perfect matching of $D-x$ by Lemma~\ref{lemWest2011}(2).
If $D$ is not covered by $M$ and $V(D) \subseteq h(\mathcal{C}_3)$, then $D$ satisfies one of the two properties listed in Fact~\ref{fact:fact1}(2)-(3), or $D$ is a compact-component of $H-X$.  
We let $y\in V(D)$ be arbitrary, and let $M_D$ be a matching of $D-y$ by Lemma~\ref{lemWest2011}(2) again. 
Then $M^{\ast}=M^{\prime}\cup (\bigcup_{D}M_D)$ is a matching of $H$.   

Let $s^{\ast}$ be the total number of components of $H-X$ that contain a splitting vertex of $H$. 
As $M$ covers every vertex of $B(X)$ from $\mathcal{Y}_3$ and at least 
$\left\lceil \frac{1}{3} (|\mathcal{Y}_1|+2|\mathcal{Y}_2|)\right \rceil$ 
vertices of $\mathcal{Y}_1 \cup \mathcal{Y}_2$, by the construction of  $M^{\ast}$, we have
\begin{eqnarray*}
    \df(H) = \odd(H-X)-|X|  
        & \le & q+s^{\ast}+|\mathcal{Y}_0\cup\mathcal{Y}_1 \cup \mathcal{Y}_2|-\lceil \frac{1}{3} (|\mathcal{Y}_1|+2|\mathcal{Y}_2|) \rceil \\
        & \le & q+s^{\ast}+|\mathcal{Y}_0|+\frac{2}{3}|\mathcal{Y}_1|+\frac{1}{3}|\mathcal{Y}_2|\\
        & \le & q+s^{\ast} +\frac{1}{3}|T_1|, \quad \text{by Fact~\ref{fact:fact1}(2)-(3) and the definition of $\mathcal{Y}_i$.}
\end{eqnarray*}

We now analysis how vertices in $h(\mathcal{C}_3)$ are covered by $M^{\ast}$. For any $x\in h(\mathcal{C}_3)$, there are three cases as follows:
\begin{enumerate}[(a)]
	\item 
    $x$ is not covered by $M^{\ast}$. 
    As $|M^{\ast}|=\frac{1}{2}(n(H)-\df(H))$ and $M$ covers $x$ if it is contained in a component of $H-X$ that contains a splitting vertex of $H$, by the upper bound on $\df(H)$ above, there are at most $ q+\frac{1}{3}|T_1|$ vertices of $h(\mathcal{C}_3)$ that are not covered by $M^{\ast}$.
	\item 
    $x$ is covered by $M^{\ast}$ but is matched to a vertex $x^{\prime}$ that is a splitting vertex of $H$.
    As $p=\sum_{v\in T_2}(c(v)-2)$, by the construction of $H$ from $G_3$, there are at most $p+\sum_{k\ge 2}\frac{1}{3}(2k+1)c_{2k+1}+\frac{2}{3}$ such vertices $x$ of $h(\mathcal{C}_3)$.
	\item 
    $x$ is covered by $M^{\ast}$ and is matched to another vertex $x^{\prime}\in  h(\mathcal{C}_3)$.
    Assume there are in total $2m$ vertices of $h(\mathcal{C}_3)$ that are covered by $M^{\ast}$ in this way.
\end{enumerate}
We can now conclude the following claim, and complete our proof in Stage 1. 
\begin{claim}\label{claim:c3-size}
    $c_3\le 2m+p+\sum_{k\ge 2} \left(\frac{2}{3}k+\frac{1}{3}\right)c_{2k+1}+\frac{1}{3}|T_1|+q+\frac{2}{3}$.
\end{claim}

\smallskip 

{\bf \noindent Stage 2: Establish a lower bound on $|S|$.}

\smallskip 
In the rest of the proof, we show that the size of $S$ is relatively large, related to a sum involving the $c_{2k+1}$ coefficients, $|T_1|$, $p$ and $q$. 
Building on this, we identify a cutset $S^{\prime}$ in $G$ that includes $S$, a small subset of vertices from $T$, and some vertices from $U$. 
We then show that $|S^{\prime}|/c(G-S^{\prime})<\frac{3}{2}$, leading to a contradiction. 
For easier description, we call the following the {\bf dist-condition}:
\begin{equation*}
    \text{The distance in $G$ between any two distinct vertices of degree 3 is at least three}. 
\end{equation*}

For $v\in V(G)$, the face of $G-v$ that contains all the vertices from $N_G(v)$ is called the \emph{induced face} of $v$ in $G$ and is denoted by $F_G(v)$.
As $G$ is a plane triangulation, $F_G(v)$ is a cycle on $N_G(v)$.
For each $v\in T$, as $N_G(v)\cap U$ is an independent set by Lemma~\ref{lemKanno2019}(2)-(3), we have $e_G(v,S)\ge e_G(v,U)=c^{\prime}(v)$. 
Denote by $f_G(v)$ the length of $F_G(v)$.
Thus
\begin{equation}\label{eq:f_G(v)}
    \text{$f_G(v)=c^{\prime}(v)+e_G(v,S)\ge 2c^{\prime}(v)\ge 2c(v)$\quad for each $v\in T$.}
\end{equation}

We first contract each component from $\mathcal{C}\setminus \mathcal{C}_1$ to a single vertex. 
Let $\mathcal{C}_1(v)=\{D\in \mathcal{C}_1: e_G(D,v)=1\}$. 
For each $v\in T$, we contract all the components in $\mathcal{C}_1(v)$ to the vertex $v$ and still call $v$ the resulting vertex. 
Denote by $G^{\prime}$ the resulting multiple plane graph.  
By Lemma~\ref{lemKanno2019}(3)-(4), multiple edges of $G^{\prime}$ exists only between $S$ and the vertices obtained by contracting components of $\mathcal{C}$. 
We remove  multiple edges in $G'$ to get a simple graph.   
In particular, for each $v\in T$, if there are multiple edges $e_1,e_2$ joining $u$ and $s$ for  $s,u\in N_{G^{\prime}}(v)$, where one of $e_1$ and $e_2$ is contained in $F_{G^{\prime}}(v)$ and the other is not, we keep exactly one edge between $s$ and $u$ that is contained in $F_{G^{\prime}}(v)$.   
Denote by $G^{\prime\prime}$ the resulting simple plane graph.  
Note that $G^{\prime\prime}$ may no longer be a triangulation as some multiple edges may have been removed.
Since $F_G(v)$ is a cycle, by our way of constructing $G^{\prime\prime}$, we have $F_{G^{\prime\prime}}(v)$ is still a cycle. 
Now if there are two nonadjacent $S$-vertices in $G^{\prime\prime}$ and we can add an edge joining them without violating the planarity, we add such an edge.  
In particular, we add such edges in the following way: if $s_1s_2\ldots s_ks_1$ for some integer  $k\ge 4$ is a facial walk in $G^{\prime\prime}$ and there exists some $i\in [1,k]$ such that $s_i, s_{i+2}$ are distinct $S$-vertices with $s_is_{i+2} \not\in E(G^{\prime\prime})$, we add the edge $s_is_{i+2}$ to $G^{\prime\prime}$ within the face of $G^{\prime\prime}$ bounded by $s_1s_2\ldots s_ks_1$. 
Denote by $G^{\ast}$ the resulting plane graph. 

For each $v\in T$, by the construction of $G^{\ast}$, we have $F_{G^{\ast}}(v)$ is a cycle on $N_{G^{\ast}}(v)$.
We claim that $d_{G^{\ast}}(v)\ge 3$ for each $v\in T$.
This is clear by \eqref{eq:f_G(v)} if $c(v)=c^{\prime}(v)$ or $c^{\prime}(v)\ge 3$.
Thus we have $c^{\prime}(v)\le 2$ and $c^{\prime}(v)-c(v)\ge 1$.
Since $G$ is 3-connected, $v$ together with all vertices from components in $\mathcal{C}_1(v)$ have in $G$ at least three distinct neighbors from $V(G)\setminus ((\bigcup_{D\in \mathcal{C}_1(v)}V(D))\cup T)$. 
Thus by the construction of $G^{\ast}$, we have $d_{G^{\ast}}(v)\ge 3$.


Let $T_{11}=\{v\in T_1: d_{G^{\ast}}(v)=3\}$ and $v\in T_{11}$. 
Then $c(v)=1$.
Since $d_{G^{\ast}}(v)\ge c(v)+e_G(v,S)$ and $e_G(v,S)\ge c(v)$, it then follows that $e_G(v,S)=2$. 
Thus we have that either $c(v)=c^{\prime}(v)=1$ ($d_G(v)=3$) or $c^{\prime}(v)-c(v)=1$ ($d_G(v)=4$). 
When $d_G(v)=3$, $F_G(v)$ and so $F_{G^{\ast}}(v)$ is a triangle. Thus $F_{G^{\ast}}(v)$ contains a unique edge $s_1s_2$ with $s_1,s_2\in S$. 
When $d_G(v)=4$ and $c^{\prime}(v)-c(v)=1$, for the component $D\in \mathcal{C}_1$ such that $e_G(v,D)=1$, we have $e_G(D,S)=2$ and $N_G(D)\cap S=N_G(v)\cap S$.
Let $N_G(D)\cap S=\{s_1,s_2\}$.
Then as $G$ is a plane triangulation, we have $s_1s_2\in E(G)\cap E(G^{\ast})$. 
Thus
\begin{equation}\label{eq:F_G*(v)}
    \text{$F_{G^{\ast}}(v)$ is a triangle and contains a unique edge $s_1s_2\in E(G)$ with $s_1,s_2\in S$.}
\end{equation}
We call triangles like $vs_1s_2v$  of $G^{\ast}$ for $v\in T_{11}$ \emph{bad triangles}. 
For $v\in T_{11}$, we denote by $F_{G^{\ast},S}(v)$ the unique edge contained in $F_{G^{\ast}}(v)$ such that the two endvertices of the edge are contained in $S$. 
We claim the following.

\begin{claim}\label{claim:cutset-bad-triangle}
	Let $v\in T_{11}$ and $F_{G^{\ast},S}(v)=s_1s_2$.   
    If $d_G(v)=4$, then $\{v,s_1,s_2\}$ is a cutset of $G$ such that one of the components of $G-\{v,s_1, s_2\}$ is trivial, and the trivial component consisting of the vertex $y^{\ast}\in U_1$ for which $e_G(v,y^{\ast})=1$.  
\end{claim}

\proof 
Let $D\in \mathcal{C}\setminus \mathcal{C}_1$ and $y\in V(D)$ such that $e_G(v, D)=e_G(v,y)=1$. 
If $d_G(v)=4$, then we further let $D^{\ast}\in \CC_1$ and $y^{\ast}\in V(D^{\ast})$ such that $e_G(v,D^{\ast})=e_G(v,y^{\ast})=1$. 
When $d_G(v)=4$, as $G$ is a plane triangulation with $yy^{\ast} \not\in E(G)$, the induced face $F_G(v)$ of $v$ is $ys_1y^{\ast}s_2y$.  
As $s_1s_2\in E(G)$ by~\eqref{eq:F_G*(v)}, it follows that $\{v,s_1,s_2\}$ is a cutset of $G$ with $V(D)$ and $V(D^{\ast})$ respectively contained in two distinct components of $G-\{v,s_1,s_2\}$. 

Assume otherwise that $G-\{v,s_1,s_2\}$ has two nontrtivial components $G_1$ and $G_2$. 
As $G$ is 3-connected, each of the component $G_i$ has a vertex that is adjacent in $G$ to $v$.  
Then $c(G-N_G(v)) \ge 3$, where one of the component contains vertices of $G_1-N_G(v)$, one contains vertices of $G_2-N_G(v)$, and one is trivial consisting of $v$. 
As $|N_G(v)|=d_G(v)=4$, we get a contradiction to $\tau(G)\ge \frac{3}{2}$.  

Thus $G-\{v,s_1,s_2\}$ has at most one nontrivial component. 
As $|V(D)| \ge 3$ by Lemma~\ref{lemKanno2019}(3)-(4), it follows that the component of $G-\{v,s_1, s_2\}$  that contains a vertex of $D^{\ast}$ is trivial.  
Thus $D^{\ast}$ is a trivial component of $G-\{v,s_1, s_2\}$ and so $V(D^{\ast})=\{y^{\ast}\}$. 
\qed

Let $v\in T_{11}$ and $F_{G^{\ast},S}(v)=s_1s_2$. 
By Claim~\ref{claim:cutset-bad-triangle}, if $d_G(v)=4$, then $v$ has a neighbor $w\in U_1$ for which $d_G(w)=3$ and $ws_1s_2w$ is a facial triangle of $G$.   
By the dist-condition, there is at most one vertex of degree 3 that is adjacent in $G$ to both $s_1$ and $s_2$.  
Thus, when $d_G(v)=3$, we call $v$ the \emph{3-neighbor} associated with $s_1s_2$; 
and when  $d_G(v)=4$, we call $w$ the \emph{3-neighbor} associated with $s_1s_2$.

\begin{claim}\label{claim:F_{G*,S}(v)}
	Any two bad triangles are vertex-disjoint in $G^{\ast}$.
\end{claim}

\proof 
Let distinct $u, v\in T_{11}$.  
When $d_G(u)=4$, we let $w_1 \in U_1$ be the 3-neighbor associated with $F_{G^{\ast},S}(u)$; and when $d_G(v)=4$, we let $w_2 \in U_1$ be the 3-neighbor associated with $F_{G^{\ast},S}(u)$. 
Since every vertex of $U_1$ is adjacent in $G$ to exactly one vertex from $T$, we have $w_1\ne w_2$. 
Furthermore, we have $u \ne w_2$ and $v \ne w_1$.  
Thus the 3-neighbor associated with  $F_{G^{\ast},S}(u)$ is distinct from that associated with $F_{G^{\ast},S}(v)$. 
By the dist-condition, each of $F_{G^{\ast},S}(u)$ and $F_{G^{\ast},S}(v)$ has a unique 3-neighbor in $G$. 
As a consequence, $F_{G^{\ast},S}(u)$ and $F_{G^{\ast},S}(v)$ are vertex-disjoint in $G$ and so in $G^{\ast}$. 
Therefore any two bad triangles are vertex-disjoint in $G^{\ast}$.
\qed

%

\begin{claim}\label{claim:3-v-in-T11b}
    For any $v\in T_{11}$, let  $F_{G^{\ast},S}(v) =s_1s_2$, and let $D\in \mathcal{C}\setminus \mathcal{C}_1$ and $y\in V(D)$ such that $e_G(v,D)=e_G(v,y)=1$. 
    If $s_1s_2$ is contained in a facial triangle $us_1s_2u$ of $G$ with $u\in V(D)$, then there exists $s^{\ast}\in S\setminus \{s_1,s_2\}$ such that $us^{\ast}s^{\prime}u$ is a facial triangle in $G$ with $s^{\ast}s^{\prime\prime}\not\in E(G)$, where $\{s^{\prime},s^{\prime\prime}\}= \{s_1,s_2\}$.  
\end{claim}

\begin{figure}[!htb]
	\begin{center}

	\begin{tikzpicture}
		
		{\tikzstyle{every node}=[draw,circle,fill=black,minimum size=4pt,
			inner sep=2pt]
			
		\draw[black, thick] (1,2.5) ellipse [x radius=2cm, y radius=1cm];
			
		\node[label=right:{$u$}] (u) at (2.9,2.2) {};
		\node[label=left:{$z_q$}] (zq) at (2.7,1.45) {};
		\node[label=left:{$z_{q+1}$}] (zq1) at (2.65,1) {};
		\node[label=below left:{$s_2$}] (s2) at (-1,0) {};
		\node[label=below right:{$s_1$}] (s1) at (3,0) {};
		\node[label=below:{$v$}] (v) at (1,0.5) {};
		\node[fill=black, inner sep=2pt, label=above:{$y$}] (y) at (1,1.5) {};
		
		\draw (y) -- (s2) -- (v) -- (s1) -- (u);
		\draw (y) -- (v);
		\draw (zq) -- (zq1);
		\draw (zq) -- (u);
		\draw (zq) -- (s1);
		\draw [red, thick](y) -- (s1);
		\draw [red, thick](u) -- (s1);
		\draw (y) -- (s2);
		\draw (s1) -- (s2);
		\draw (u) .. controls (6,-0.5) and (3,-2) .. (s2);
		\draw [thick,red] (y) .. controls (1.7,1.5) and (2.5,1.66) .. (u);
		
	}

		\node[draw=none, fill=none, text=black] at (1,2.5) {\Large $D$};
	\end{tikzpicture}

\end{center}
\vspace{-0.5cm}
\caption{An illustration of the embedding of $G$ around the vertex $v$, where $\Ints(vs_1s_2v)$ is either a null graph or is a single vertex graph, and the boundary drawn in red represents $W$.}
\label{f1}
\end{figure}
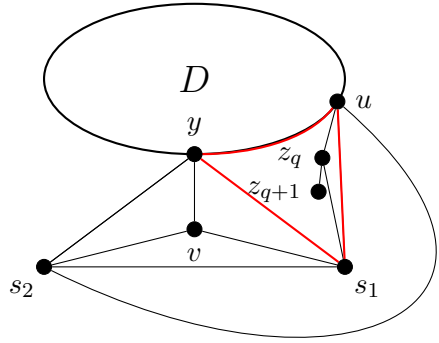

\proof  
Since $G$ is plane triangulation and $d_G(v) \le 4$ with $N_G(v)\cap V(D)=\{y\} $ by Lemma~\ref{lemKanno2019}(3), we have $vys_1v$ and $vys_2v$ are both facial triangles in $G$.  
Let $u\in V(D)$ such that $us_1s_2u$ is a facial triangle in $G$. Without loss of generality, we assume that the embedding of $G$ is as depicted in Figure~\ref{f1}. 
Let $W$ be the walk bounded by $ys_1,s_1u$ and the  outer facial walk of $D$ between $y$ and $u$ (the curves drawn in red). 
If $\Ints(W)$ is a null graph or the outer facial walk of $\Ints(W)$ consists of $S$-vertices only,  then $s_1$ is adjacent in $G$ to only one component that is $D$ from $\mathcal{C}\setminus \mathcal{C}_1$ ($\Ints(vs_1s_2v)$ is either a null graph or is a single vertex graph by Claim~\ref{claim:cutset-bad-triangle}).  
Let $S^{\ast}=S\setminus \{s_1\}$ and $T^{\ast}=T\setminus \{v\}$.  
Then we get 
\begin{eqnarray*}
    \delta(S^{\ast},T^{\ast})
        & \le & 2|S|-2-2|T|+2+\left(\sum_{w\in T}d_{G-S}(w)-1-e_G(v,U_1)\right) \\
        && -\left(c(S,T)-1-e_G(v,U_1)\right)\\
        & = & 2|S|-2|T|+\sum_{w\in T}d_{G-S}(w)-c(S,T)\\
        & = & \delta(S,T),
\end{eqnarray*}
where $e_G(v,U_1)=1$ if and only if $d_G(v)=4$. 
Thus $(S^{\ast},T^{\ast})$ is also a barrier of $G$  with $|T^{\ast}|<|T|$, a contradiction to the choice of $(S,T)$.

Thus we assume that the outer facial walk of $\Ints(W)$ contains a vertex from $T\cup U$.   
Let $W^{\ast}=z_1z_2\ldots z_kz_1$ be the outer facial walk of $\Ints(W)$, where $k\ge 3$ is an integer (here we have $k\ge 3$, as one of $z_1, \ldots, z_k$, say $z_j$, is a vertex from $T\cup U$,  which is adjacent in $G$ only  to $s_1$ from vertices of $W$ and $z_js_1$ is contained in two facial triangles in $G$).  
We let $p,q\in [1,k]$ with $p \le q$ such that $q-p$ is maximum and each vertex of $z_{p},z_{p+1},\ldots,z_q$ is an $S$-vertex and adjacent in $G$ to $D$.   
Since $z_{p-1},z_{q+1}$ are vertices from $T\cup U$, they are nonadjacent in $G$ to any of vertex  of $D$.  
As $G$ is a plane triangulation, one of $z_p$ and $z_q$ is adjacent in $G$ to $y$ and the other is adjacent in $G$ to $u$. 
Without loss of generality, we assume that $z_qu \in E(G)$. 
As $G$ is a plane triangulation, $z_qu\in E(G)$, and $z_{q+1}u\not\in E(G)$, it follows that $z_{q}s_1\in E(G)$. Thus $uz_{q}s_1u$ is a facial triangle in $G$.  
Now, letting $s^{\ast}=z_q$ gives the desired $S$-vertex.  
\qed 

\begin{claim}\label{claim:3-v-in-T11c}
    Let  $v\in T_{11}$ and  $F_{G^{\ast},S}(v) =s_1s_2$. 
    Then $s_1s_2$ is contained in exactly two distinct facial triangles of $G^{\ast}$, where one of them is $vs_1s_2v$.  
\end{claim}

\proof 
For a component $R\in \mathcal{C}\setminus \mathcal{C}_1$, we let $v_R$ denote the vertex of $G^{\ast}$ obtained from contracting $R$.
Let $D\in \mathcal{C}\setminus \mathcal{C}_1$ and $y\in V(D)$ such that $e_G(v, D)=e_G(v,y)=1$. 
If $d_G(v)=4$, then we further let $D^{\ast}\in \mathcal{C}_1$ and $V(D^{\ast})=\{y^{\ast}\}$ such that $e_G(v, D^{\ast})=e_G(v,y^{\ast})=1$. 

When $d_G(v)=3$, it is clear that $vs_1s_2v$ is a facial triangle in $G$ and so in $G^{\ast}$.  
When $d_G(v)=4$, as $y^{\ast}s_1s_2y^{\ast}$ is a facial triangle in $G$ by Claim~\ref{claim:cutset-bad-triangle}, and the edge $y^{\ast}v$ is contracted into $v$, it follows that $vs_1s_2v$ is a facial triangle in $G^{\ast}$. 
We next show that $s_1s_2$ is contained in a second facial triangle of $G^{\ast}$. 

As $s_1s_2\in E(G)$ by~\eqref{eq:F_G*(v)}, there exists $u\in V(G)$ such that $us_1s_2u$ is a facial triangle of $G$ with $u\ne v$ when $d_G(v)=3$ and  with $u\ne y^{\ast}$ when $d_G(v)=4$.

Consider first that $u\not\in V(D)$. 
If $u\in S\cup T$, then $us_1s_2u$ is still a facial triangle of $G^{\ast}$.  
Thus we assume $u\in U$, and let $u\in V(Q)$ for some $Q\in \mathcal{C}$.
If $Q\not\in \mathcal{C}_1$, then $v_Qs_1s_2v_Q$ is a facial triangle of $G^{\ast}$. 
Thus we assume that $Q\in \mathcal{C}_1$, and let $z\in T$ such that $e_G(z,Q)=1$. 
As $u\ne y^{\ast}$, we have $z\ne v$. Then  $zs_1s_2z$ is a facial triangle of $G^{\ast}$.  

We consider now that $u\in V(D)$. 
By Claim~\ref{claim:3-v-in-T11b}, there exists $s^{\ast}\in S\setminus \{s_1,s_2\}$, and say $s_1$, such that $ us^{\ast}s_1u$ is another facial triangle of $G$ containing $us_1$ and $s^{\ast}s_2\not\in E(G)$.  
Thus there are at least two edges between $v_D$ and each of $s_1$ and $s_2$ 
in $G^{\prime}$, which recall was obtained from $G$ by contracting components from $\mathcal{C}$. When we got $G^{\prime\prime}$ from $G^{\prime}$, we deleted the multiple edges $v_Ds_1$ and $v_Ds_2$ such that they are not contained in the facial triangles of $G^{\prime}$ containing the vertex $v$.  
Now as $ us^{\ast}s_1u$ is a facial triangle in $G$ and $us_1s_2u$ is another facial triangle in $G$, in $G^{\prime\prime}$, $v_Ds^{\ast}s_1s_2$ is a section of a facial walk $W^{\prime}$ with  $s^{\ast}s_2\not\in E(G)$.  
Thus, by the construction of $G^{\ast}$ from $G^{\prime\prime}$, there exists $s^{\ast\ast} \in S\setminus \{s_1,s_2\}$, where $s^{\ast\ast}$ is either $s^{\ast}$ or the second vertex immediately following $s_2$ in the direction from $v_D$ to $s_2$ passing $s^{\ast}$ on $W^{\prime}$, such that $s^{\ast\ast}s_1s_2s^{\ast\ast}$ is a facial triangle of $G^{\ast}$. 
\qed 

By Claim~\ref{claim:3-v-in-T11c}, for each $v\in T_{11}$, we let $w$ be the other vertex of $G^{\ast}$ such that $F_{G^{\ast},S}(v)$ is contained in the two facial triangles respectively containing $v$ and $w$. 
We call $w$ the \emph{co-triangle vertex} of $v$, and call $(v,w)$ a \emph{co-triangle pair}.  
By Claim~\ref{claim:F_{G*,S}(v)}, for each vertex $v\in T_{11}$ with $F_{G^{\ast},S}(v)=s_1s_2$ and $w$ as its co-triangle vertex, there are two distinct edges $ws_1,ws_2$ of $G^{\ast}$ that are uniquely associated with $v$. 
For $D\in \mathcal{C}$, we let $h_2(D)$ be the vertex of $G^{\ast}$ that is obtained by contracting $D$. 
For $\mathcal{Z}\subseteq \mathcal{C}$, we also let $h_2(\mathcal{Z}) =\{h_2(D): D\in \mathcal{Z}\}$.

We partition $T_{11}$ into three subsets as follows.
Let 
\begin{equation*}
\begin{split}
    T_{11}^1 & = \{v\in T_{11}: \text{the co-triangle vertex of $v$ belongs to $S$}\}, \\
    T_{11}^2 & = \{v\in T_{11}: \text{the co-triangle vertex of $v$  belongs to $T\setminus T_{11}$}\},\quad \text{and}\\
    T_{11}^3 & = \{v\in T_{11}: \text{the co\text-triangle vertex of $v$ belongs to $h_2(\mathcal{C}\setminus \mathcal{C}_1)$}\}.
\end{split}
\end{equation*}
By Claim~\ref{claim:F_{G*,S}(v)}, the co-triangle vertex of $v$ is not any vertex from $T_{11}$. 
Thus, indeed, we have $T_{11}^1\cup T_{11}^2\cup T_{11}^3=T_{11}$.
Corresponding to $T_{11}^1$, $T_{11}^2$ and $T_{11}^3$, we define the following sets: 
\begin{equation*}
\begin{split}
    S^{\ast} & = \{w\in S: \text{$v\in T^1_{11}$ and $(v,w)$ is a co-triangle pair} \},\\
    T^{\ast} & = \{w\in T\setminus T_{11}: \text{$v\in T^2_{11}$ and $(v,w)$ is a co-triangle pair}\},\quad \text{and}\\
    U^{\ast} & = \{w\in h_2(\mathcal{C}\setminus \mathcal{C}_1): \text{$v\in T^3_{11}$ and $(v,w)$ is a co-triangle pair}\}.
\end{split}
\end{equation*}
For a set $Z$  and $w\in Z$, we let  
$\alpha(w)=|\{ v\in T_{11}: \text{$(v,w)$ is a co-triangle pair}\}|$. 
Then we have 
$$
    \sum_{w\in S^{\ast}}\alpha(w)=|T_{11}^1|, \quad 
    \sum_{w\in T^{\ast}}\alpha(w)=|T_{11}^2|, \quad \text{and} \quad 
    \sum_{w\in U^{\ast}}\alpha(w)=|T_{11}^3|.
$$
We define $T^{\ast}$ as follows: 
$$
    T_0^{\ast}=T_0\cap T^{\ast}, \quad 
    T_1^{\ast}=T_1\cap T^{\ast}, \quad \text{and}\quad 
    T_2^{\ast}=T_2\cap T^{\ast}. 
$$
For notation simplicity, for any set $Z$, we let 
$$
    m(Z)=\sum_{w\in Z}\alpha(w).
$$
It is clear that $m(Z)\ge |Z|$ if $Z\subseteq S^{\ast} \cup T^{\ast}\cup U^{\ast}$. 

By Claim~\ref{claim:componentH-X}, for any compact-component $D$ of $H-X$, the subgraph $h_1^{-1}(D)$ of $G$ has at least two auxiliary walks $W_1$ and $W_2$, and $V(W_1),V(W_2) \subseteq S\cup U_1$.  
As $G$ is a plane triangulation, and every vertex of $U_1$ has exactly one neighbor from $T$ and all its other neighbors are from $S$, it follows that the subgraph of $G^{\ast}$ that is corresponding to $W_i$  is still a closed walk. 
For notation simplicity, we also denote by $h_2(W_i)$ the subgraph of $G^{\ast}$ that is corresponding to $W_i$.

Let $W$ be an auxiliary walk of a compact-component $D$ of $H-X$.  
If $V(W) \subseteq S$ or every $u\in V(W)\cap U_1$ satisfying  $|N_G(u)\cap S|=2$, then by our construction of $G^{\ast}$, $V(h_2(W))$ consists of $S$-vertices only.  
By Claim~\ref{claim:componentH-X} and the construction of $G^{\ast}$, each $S$-vertex of $h_2(W)$ is contained in a cycle of $h_2(W)$. 
By Claim~\ref{claim:F_{G*,S}(v)}, we have $h_2(W)$ contains an edge $s_1s_2$ with $s_1,s_2\in S$ such that  $s_1s_2$ is not contained in any bad triangles.  
Let $x\in V(h_1^{-1}(D))$ such that $xs_1s_2x$ is a facial triangle in $G$ (so there is correspondingly a facial triangle containing $s_1s_2$ and a vertex of $D$ in $G^{\ast}$).  
If $x\in T$, then we call $W$ a \emph{$T$-contributive} auxiliary walk. 
Otherwise, $x\in U_3$,  and we  call $W$ a \emph{$U$-contributive} auxiliary walk. 
We let $\mathcal{W}_1$ be the set of auxiliary walks of compact-components $D$ that either have a vertex $u\in U_1$ satisfying $|N_G(u) \cap S| \ge 3$, or are $T$-contributive; 
and let $\mathcal{W}_2$ be the set of $U$-contributive auxiliary walks. 
We have 
$$
    |\mathcal{W}_1|+|\mathcal{W}_2| \ge 2q
$$
by Claim~\ref{claim:componentH-X}.

\begin{claim}\label{claim:degree-of-T0*}
    For any $w\in T_0^{\ast}$, we have $d_{G^{\ast}}(w)\ge 4$.
\end{claim}

\proof 
Let $v\in T_{11}$ be a vertex whose co-triangle vertex is $w$, 
$D$ be the  component from $\mathcal{C}\setminus \mathcal{C}_1$ containing a vertex $u$ adjacent in $G$ to $v$, and $F_{G^{\ast},S}(v)=s_1s_2$.
Suppose to the contrary that $d_{G^{\ast}}(w)=3$. 
By Claim~\ref{claim:cutset-bad-triangle}, there exists a unique 3-neighbor associated with $F_{G^{\ast},S}(v)$.  Thus $d_G(w)=4$. 
Furthermore, by the construction of $G^{\ast}$, we have $c^{\prime}(w)=1$ and the odd component $D_w$ from $\mathcal{C}_1(w)$ satisfies $N_G(D_w)\cap S \subseteq N_G(w)\cap S$. 
Let $N_G(w)\cap (S\setminus \{s_1,s_2\})=\{s_3\}$.                                        
Since $|S|\ge2$, we have $|T|\ge 3$ by~\eqref{eq:T-size}.
Therefore $c(G-\{u,s_1,s_2,s_3\})\ge 3$, where one component of $G-\{u,s_1,s_2,s_3\}$ contains $v$, one contains $w$ and the rest contains vertices from $T\setminus \{v,w\}$. This gives a contradiction to $\tau(G)\ge \frac{3}{2}$.
\qed

\begin{claim}\label{claim:edges-in-S}
We have 
$e(G^{\ast}[S]) \ge |T_{11}|+|T_{11}^1|+\frac{1}{4}m(T_1^{\ast})+\frac{1}{4}|T_1\setminus T_{11}|+\frac{1}{2}m(T_0^{\ast})+|T_0|-|T_0^{\ast}|-\frac{1}{2}\sum_{w\in T_1^{\ast}, \alpha(w) \ge 2} 1$. 
\end{claim}

\proof  
We estimate the number of edges within $S$ in $G^{\ast}$ by considering the following possibilities: 
\begin{itemize}
	\item edges corresponding to $F_{G^{\ast},S}(v)$ for $v\in T_{11}$:  
        there are $|T_{11}|$ such edges by~\eqref{eq:F_G*(v)}. 
	\item edges incident with vertices of  $S^{\ast}$: 
        for $w\in S^{\ast}$, let $(v,w)$ be a co-triangle pair, and let $F_{G^{\ast},S}(v)=s_1s_2$.  
        Then by Claim~\ref{claim:F_{G*,S}(v)}, we have $ws_1, ws_2 \neq F_{G^{\ast},S}(u)$ for any $u\in T_{11}$, but $ws_1, ws_2 \in E(G^{\ast}[S])$.  
        There are $\alpha(w)$ pairs of such edges. 
	\item edges within the neighborhood of vertices of $T_0^{\ast}$:  
        for $w\in T_0^{\ast}$, $F_{G^{\ast}}(w)$ has exactly $d_{G^{\ast}}(w)$ edges from $E(G^{\ast}[S])$.
        By Claim~\ref{claim:degree-of-T0*}, we have $d_{G^{\ast}}(w)\ge 4$. 
        Thus, we have $d_{G^{\ast}}(w) \ge 2\alpha(w)+2$ if $\alpha(w)=1$. 
        If $\alpha(w)\ge 2$, then we have $d_{G^*}(w) \ge 2\alpha(w)$ by Claim~\ref{claim:F_{G*,S}(v)}. 
        Again, by Claim~\ref{claim:F_{G*,S}(v)}, at least $\left \lceil \frac{d_{G^{\ast}}(w)}{2}  \right \rceil $ edges of $F_{G^{\ast}}(w)$ are not $F_{G^{\ast},S}(u)$ for any $u\in T_{11}$. 
	 \item edges within the neighborhood of  vertices of  $T_1^{\ast}$: 
        for  $w\in T_1^{\ast}$, let $(v,w)$ be a co-triangle pair, and let $F_{G^{\ast},S}(v)=s_1s_2$.   
        Since $T_1^{\ast}\cap T_{11} =\emptyset$, we have  $|N_G(w)\cap S|\ge 3$. 
        If $\alpha(w)=1$, let $s_3\in N_G(w)\cap (S\setminus \{s_1,s_2\})$ such that $s_1s_3\in E(G)$. 
        Then $s_1s_3\neq F_{G^{\ast},S}(u)$ for any $u\in T_{11}$ by Claim~\ref{claim:F_{G*,S}(v)}. 
        If $\alpha(w) \ge 2$, then $|N_G(w)\cap S|\ge 2\alpha(w)$ by Claim~\ref{claim:F_{G*,S}(v)}. Thus, when $\alpha(w)=1$, at least one edge of $E(F_{G^{\ast}}(w))\cap E(G^{\ast}[S])$ is not $F_{G^{\ast},S}(u)$ for any $u\in T_{11}$.
        When $\alpha(w) \ge 2$, at least $\alpha(w)-1$ edges of $E(F_{G^{\ast}}(w))\cap E(G^{\ast}[S])$ are not $F_{G^{\ast},S}(u)$ for any $u\in T_{11}$. 
	 \item edges within the neighborhood of vertices of $T_1\setminus (T_{11} \cup T_1^{\ast})$:  
        let $v\in T_1\setminus (T_{11} \cup T_1^{\ast})$. 
        Note that $d_{G^{\ast}}(v)\ge4$.
        Then $F_{G^{\ast}}(v)$ has exactly $d_{G^{\ast}}(v)-2$ edges from $E(G^{\ast}[S])$. 
        As $v\notin T_1^{\ast}$, for each edge $e\in E(F_{G^{\ast}}(v))\cap E(G^{\ast}[S])$, we have $e\neq F_{G^{\ast},S}(u)$ for any $u\in T_{11}$. 
    \item edges within the neighborhood of  vertices of $T_0\setminus T_0^{\ast}$:  
        let $v\in T_0\setminus T_0^{\ast}$.  
        Then $F_{G^{\ast}}(v)$ has exactly $d_{G^{\ast}}(v)$ edges from $E(G^{\ast}[S])$ which are not $F_{G^{\ast},S}(u)$ for any $u\in T_{11}$.
\end{itemize}
As an edge from $E(G^{\ast}[S])$ can be counted twice in all the cases above except the first one,  thus $e(G^{\ast}[S])$ is at least 
\begin{eqnarray*}
    && \sum_{v\in T_{11}}|\{ F_{G^{\ast},S}(v)\}|+\frac{1}{2}\sum_{w\in S^{\ast}}2\alpha(w)+\frac{1}{2}\sum_{w\in T_1^{\ast}, \alpha(w)=1}1+ \frac{1}{2}\sum_{w\in T_1^{\ast}, \alpha(w) \ge 2}(\alpha(w)-1) \\
    && +\frac{1}{2}\sum_{w\in T_0^{\ast}, \alpha(w)=1}\frac{2\alpha(w)+2}{2}+\frac{1}{2}\sum_{w\in T_0^{\ast}, \alpha(w) \ge 2} \frac{2\alpha(w)}{2}  \\ 
    && +\frac{1}{2}\sum_{v\in T_1\setminus (T_{11} \cup T^{\ast})}(d_{G^{\ast}}(v)-2)+\frac{1}{2}\sum_{v\in T_0\setminus T_0^{\ast}}d_{G^{\ast}}(v)\\
    & \ge & |T_{11}|+|T_{11}^1|+\frac{1}{2}m(T_1^{\ast})-\frac{1}{2}\sum_{w\in T_1^{\ast}, \alpha(w) \ge 2} 1+ \frac{1}{2}m(T_0^{\ast})+|T_1\setminus (T_{11} \cup T_1^{\ast})|+\frac{3}{2}|T_0\setminus T_ 0^{\ast}|\\
    & \ge & |T_{11}|+|T_{11}^1|+\frac{1}{4}m(T_1^{\ast})+|T_1\setminus T_{11}|-\frac{3}{4}|T_1^{\ast}|+\frac{1}{2}m(T_0^{\ast})+|T_0|-|T_0^{\ast}|-\frac{1}{2}\sum_{w\in T_1^{\ast}, \alpha(w) \ge 2} 1\\
    &\ge & |T_{11}|+|T_{11}^1|+\frac{1}{4}m(T_1^{\ast})+\frac{1}{4}|T_1\setminus T_{11}|+\frac{1}{2}m(T_0^{\ast})+|T_0|-|T_0^{\ast}|-\frac{1}{2}\sum_{w\in T_1^{\ast}, \alpha(w) \ge 2} 1, 
\end{eqnarray*}
where we have $|T_1\setminus T_{11}|-\frac{3}{4}|T_1^{\ast}|\ge \frac{1}{4}|T_1\setminus T_{11}|$ because $T_1^{\ast}\subseteq T_1\setminus T_{11}$. 
\qed

\begin{claim}\label{claim:S-U-edges}
We have 
$e_{G^{\ast}}(S, h_2(\mathcal{C}\setminus \mathcal{C}_1)) \ge \sum_{k\ge 1}(2k+1)c_{2k+1}+|T_{11}^3|+|\mathcal{W}_2|$. 
\end{claim}

\proof 
For any $w\in U^{\ast}$, let $v$ be a vertex from $T_{11}$ whose co-triangle vertex is $w$, and  let $F_{G^{\ast},S}(v) =s_1s_2$.  
Then as $vs_1s_2v$ is a facial triangle of $G^{\ast}$, each of $ws_1$ and $ws_2$ is contained in at most one induced face of some vertex from $T$. 
There are $\alpha(w)$ pairs of such edges. 
For any $W\in \mathcal{W}_2$, $W$ has an edge $s_1s_2$  with $s_1,s_2\in S$ for which there is $x\in h_2(\mathcal{C}_3)$  such that $xs_1s_2x$ is a facial triangle of $G^{\ast}$, and that  $x\not\in U^*$.   
Then each of $xs_1$ and $xs_2$ is contained in at most one induced face of some vertex from $T$ in $G^{\ast}$. 
Therefore, $e_{G^{\ast}}(S,h_2(\mathcal{C}\setminus \mathcal{C}_1))$ is at least 
\begin{eqnarray*}
    && \frac{1}{2}\sum_{v\in T}e_{G^{\ast}}(S\cap V(F_{G^{\ast}}(v)),h_2(\mathcal{C}\setminus \mathcal{C}_1)\cap V(F_{G^{\ast}}(v)))+\frac{1}{2}\sum_{w\in U^{\ast}}2\alpha(w)+\frac{1}{2}\sum_{W\in \mathcal{W}_2}2\\ 
    & \ge & \frac{1}{2}\sum_{v\in T}2c(v)+m(U^{\ast})+|\mathcal{W}_2|
    = \sum_{k\ge 1}(2k+1)c_{2k+1}+|T_{11}^3|+|\mathcal{W}_2|.
\end{eqnarray*}
\qed

\begin{claim}\label{claim:S-size}
We have 
$|S| \ge \sum_{k\ge 1}\left (\frac{4}{3}k-\frac{5}{6}\right )c_{2k+1}+\frac{1}{6}|T_1|+\frac{5}{6}p+\frac{1}{2}q+3$. 
\end{claim}

\proof  
We start by imposing a lower bound on $|S|$ by counting the edges in $G^{\ast}-T$.  
By Claim~\ref{claim:degree-of-T0*} and the definition of $T_1\setminus T_{11}$ and $T_2$, we have $d_{G^{\ast}}(v) \ge 4$ for any $v\in T_0^{\ast} \cup (T_1\setminus T_{11})\cup T_2$. 
Thus $G^{\ast}-T$ has at least $|T_0^{\ast}|+|T_1\setminus T_{11}|+|T_2|$ faces of length at least 4. 
Note that $f_G(v)\ge 2c(v)$ for each $v\in T$ by~\eqref{eq:f_G(v)} and $f_{G^{\ast}}(u)\ge 2c(u)+\alpha(u)$ for each $u\in T_2^{\ast}$.
Furthermore, for each $w\in T_1^{\ast}$ with $\alpha(w) \ge 2$, we have $f_{G^{\ast}}(w) \ge 2\alpha(w)+1$. 
For $W\in \mathcal{W}_1$, let $D$ be the compact-component of $H-X$ such that  $W$ is 
an auxillary walk of $D$. 
If there exists $u\in V(W)\cap U_1$ such that $|N_G(u) \cap S|\ge 3$,  we let $x\in V(h_1^{-1}(D))\cap T_2$ such that $xu\in E(G)$. 
Then by the construction of $G^{\ast}$, there exist three consecutive $S$-vertices on $F_{G^{\ast}}(x)$.  
As a consequence, we have $d_{G^{\ast}}(x) \ge 2c(x)+2$.  
We denote by $T_2^{a}$ the set of such $T_2$-vertices $x$.  
For any $x\in T_2^{a}\cap T_2^{\ast}$, $d_{G^{\ast}}(x) \ge 2c(x)+2$ if $\alpha(x) \le 1$ and $d_{G^{\ast}}(x) \ge 2c(x)+\alpha(x)+1$ if $\alpha(x) \ge 2$.  
In both cases, we have $f_{G^{\ast}}(x)-4=d_{G^{\ast}}(x)-4 \ge 2c(x)-4+\alpha(x)+1$. 
If $W$ is a $T$-contributive auxillary walk, then  $W$ has an edge $s_1s_2$ with $s_1,s_2\in S$ such that $s_1s_2$ is not contained in any bad triangles. 
Let $x\in V(h_1^{-1}(D))\cap T_2$ such that $xs_1s_2x$ is a facial triangle of $G$. 
Then $xs_1s_2x$ is also a facial triangle of $G^{\ast}$, and we have $d_{G^{\ast}}(x) \ge 5$ and $x\not\in T_2^{\ast}$.   
Denote by $T_2^{b}$ the set of such $T_2$-vertices $x$. 
Thus, we have
\begin{eqnarray*}
	&& \sum_{v\in T_2}(f_{G^{\ast}}(v)-4) \\
    & \ge & \sum_{v\in T_2}(2c(v)-4)+ \sum_{x\in T_2^{\ast} \cap T_2^{a}}(\alpha(x)+1)+m(T_2^{\ast} \setminus  T_2^{a})+2|T_2^a\setminus T_2^{\ast}|+|T_2^b| \\
	& \ge & 2p+m(T_2^{\ast})+|T_2^a|+|T_2^b|=2p+m(T_2^{\ast})+|\mathcal{W}_1|. 
\end{eqnarray*}
Therefore, 
\begin{eqnarray}
    && e_{G^{\ast}}(S,V(G^{\ast})\setminus (S\cup T)) +e(G^{\ast}[S])  \label{eq:sum-SU-S}  \\
    & \le & e(G^{\ast}-T)
    \le 3|V(G^{\ast}-T)|-6-\sum_{v\in T}(f_{G^{\ast}}(v)-3) \nonumber \\
    & \le & 3|S| +\sum_{k\ge 1}3c_{2k+1}-6-|T_0^{\ast}|-|T_1\setminus T_{11}|-|T_2|-m(T_2^{\ast})-|\mathcal{W}_1|-2p \nonumber \\
    && - \sum_{w\in T_1^{\ast},\alpha(w) \ge 2} (2\alpha(w)+1-3-1) \nonumber  \\
    & \le & 3|S| +\sum_{k\ge 1}3c_{2k+1}-6-|T_0^{\ast}|-|T_1\setminus T_{11}|-|T_2|-m(T_2^{\ast})-|\mathcal{W}_1|-2p-\sum_{w\in T_1^{\ast},\alpha(w) \ge 2}1.\nonumber 
\end{eqnarray}

Now, Claim~\ref{claim:edges-in-S}, Claim~\ref{claim:S-U-edges} and~\eqref{eq:sum-SU-S} together give
\begin{eqnarray}
    && |T_{11}|+|T_{11}^1|+\frac{1}{4}m(T_1^{\ast})+\frac{1}{4}|T_1\setminus T_{11}|+ \frac{1}{2}m(T_0^{\ast})+|T_0|-|T_0^{\ast}| \nonumber \\
    && + \sum_{k\ge 1}(2k+1)c_{2k+1}+|T_{11}^3|+|\mathcal{W}_2| \nonumber \\
    & \le & 3|S|+\sum_{k\ge 1}3c_{2k+1}-6-|T_0^{\ast}|-|T_1\setminus T_{11}|-|T_2|-m(T_2^{\ast})-|\mathcal{W}_1|-2p. \label{eq:sum-US-SS-sum}
\end{eqnarray}
Since $|\mathcal{W}_1|+|\mathcal{W}_2| \ge 2q $, inequality~\eqref{eq:sum-US-SS-sum}, together with~\eqref{eq:T2-size}, gives 
\begin{eqnarray*}
    3|S|
        & \ge & |T_{11}|+|T_{11}^1|+\frac{1}{4}m(T_1^{\ast})+\frac{1}{4}|T_1\setminus T_{11}|+\frac{1}{2}m(T_0^{\ast})+ |T_0|-|T_0^{\ast}|+\sum_{k\ge 1}(2k+1)c_{2k+1} \\
        && +|T_{11}^3|-\sum_{k\ge 1}3c_{2k+1}+|T_0^{\ast}|+|T_1\setminus T_{11}|+|T_2|+m(T_2^{\ast})+2p+2q+6 \\
        & = & \sum_{k\ge 1}\left(3k-\frac{3}{2} \right)c_{2k+1}+\frac{3}{4}|T_1|-\frac{1}{4}|T_{11}| \\
        && + \left(|T_{11}^1|+\left(\frac{1}{2}m(T_0^{\ast})+\frac{1}{4}m(T_1^{\ast})+m(T_2^{\ast}) \right)+|T_{11}^3|\right)+|T_0|+\frac{3}{2}p+2q+6. 
\end{eqnarray*}
Since $\frac{1}{2}m(T_0^{\ast})+\frac{1}{4}m(T_1^{\ast})+m(T_2^{\ast})\ge \frac{1}{4}m(T^{\ast}) =\frac{1}{4}|T^2_{11}|$, and $|T^1_{11}|+|T^2_{11}|+|T^3_{11}| =|T_{11}|$, 
we then get 
\begin{equation*}
    |S|\ge \sum_{k\ge 1}\left (k-\frac{1}{2} \right)c_{2k+1} +\frac{1}{3}\left |T_0\right | +\frac{1}{4}\left |T_1\right |+\frac{1}{2} p+\frac{2}{3}q+2.
\end{equation*}
This together with~\eqref{eq:T-size} gives
\begin{equation}\label{eq:SinT}
    |T|\ge \sum_{k\ge 1}\left(2k-\frac{1}{2} \right)c_{2k+1} +\frac{1}{3}\left |T_0\right |+\frac{1}{4}\left |T_1\right | +\frac{1}{2} p+\frac{2}{3}q+3.
\end{equation}
Now by~\eqref{eq:T2-size} and~\eqref{eq:SinT}, we get
\begin{eqnarray*}
    |T_1| & = & |T|-|T_2|-|T_0| \\
          & \ge & \sum_{k\ge 1}\left (2k-\frac{1}{2} \right )c_{2k+1}+\frac{1}{3}|T_0|+\frac{1}{4}|T_1| +\frac{1}{2} p+\frac{2}{3}q+3\\
          & \quad & -\sum_{k\ge 1}\left(k+\frac{1}{2}\right )c_{2k+1}+\frac{1}{2}|T_1|+\frac{1}{2}p-|T_0| \\
          & = & \sum_{k\ge 1}(k-1)c_{2k+1}-\frac{2}{3}|T_0|+\frac{3}{4}| T_1|+p+\frac{2}{3}q+3.
\end{eqnarray*}
Therefore
\begin{equation}\label{eq:T1-size}
    |T_1| \ge \sum_{k\ge 1}(4k-4)c_{2k+1}-\frac{8}{3}\left |T_0\right |+4p+\frac{8}{3}q+12.
\end{equation}
By the definition of $\delta(S,T)$ and \eqref{eq:T-size}, we have
\begin{eqnarray*}
    2|S| & = & 2|T|-\sum_{y\in T}d_{G-S}(y)+\sum_{k\ge 0}c_{2k+1}-2 \\
         & = & 2|T_0|+2|T_1|+2|T_2|-|T_1|-2|T_2|-p-c_1+\sum_{k\ge 1}c_{2k+1}+c_1-2 \\
         & = & 2|T_0|+|T_1|+\sum_{k\ge 1}c_{2k+1}-p-2.
\end{eqnarray*}
Thus by \eqref{eq:T1-size}, we get
\begin{eqnarray*}
    |S| & = & |T_0|+\frac{1}{2}|T_1|+\sum_{k\ge 1}\frac{1}{2}c_{2k+1}-\frac{p}{2}-1 \\
        & = & |T_0|+\frac{1}{3}|T_1|+\frac{1}{6}|T_1|+\sum_{k\ge 1}\frac{1}{2}c_{2k+1}-\frac{p}{2}-1 \\
        & \ge & \sum_{k\ge 1}\left (\frac{4}{3}k-\frac{5}{6}\right )c_{2k+1}+\frac{1}{6}|T_1|+\frac{5}{6}p+\frac{8}{9}q+3,
\end{eqnarray*}
proving the claim.
\qed

For $k\ge 2$, we let
\begin{equation*}
\begin{split}
    \mathcal{C}_{2k+1}^1 & = \{D\in \mathcal{C}_{2k+1}: |V(D)|=2k+1\},\quad 
    \mathcal{C}_{2k+1}^2=\mathcal{C}_{2k+1}\setminus \mathcal{C}_{2k+1}^1,\\
    c_{2k+1}^1 & = |\mathcal{C}_{2k+1}^1|,\quad 
    c_{2k+1}^2=|\mathcal{C}_{2k+1}^2|.
\nonumber
\end{split}
\end{equation*}
We have $|V(D)|\ge 2k+2$ for each $D\in \mathcal{C}_{2k+2}^2$ by Lemma~\ref{lemKanno2019}(3)-(4).

\begin{claim}\label{claim:C_{2k+1}^1}
    Let $k\ge 2$ and $D\in \mathcal{C}_{2k+1}^1$. 
    Then there exists $W\subseteq V(D)$ with $|W| =\left \lfloor \tfrac{4k+2}{3}  \right \rfloor$ such that $V(D)-W$ is a set of $\left \lceil \tfrac{2k+1}{3}  \right \rceil $ isolated vertices.
\end{claim}

\proof
By Lemma~\ref{lemKanno2019}(3)-(4), each vertex of $D$ is adjacent in $G$ to a distinct vertex of $T$.
Thus $D$ is outer-planar. 
As a consequence, $D$ is 3-colorable and so a largest color class of $D$ has size at least $\frac{2k+1}{3}$. 
So we can choose $W\subseteq V(D)$ with $|W|=\left \lfloor \tfrac{4k+2}{3} \right \rfloor$ such that $V(D)-W$ is a set of isolated vertices.
\qed

For each $w_1w_2\in M^{\ast}$ with $w_1,w_2\in h(\mathcal{C}_3)$, by the definition of $H$, for each $i\in [1,2]$, there exists a vertex $x_i\in V (h^{-1}(w_i))$ such that $x_1$ and $x_2$ are both adjacent in $G$ to a vertex $y\in T$ with $c(y)=2$. 
We call $x_i$ a \emph{representative vertex} of $h^{-1}(w_i)$. 
As $M^{\ast}$ is a matching of $H$, each component from $\mathcal{C}_3$ either has no representative vertex or has a unique representative vertex.

Let $D\in \mathcal{C}_3$. If $D$ has a representative vertex, say $x$, let $S_D\subseteq V(D)\setminus \{x\}$ be the set of two vertices such that $e_G(D-S_D,T)=e_G(x,T)$. Otherwise, let $S_D\subseteq V(D)$ be a set of two vertices such that $e_G(D-S_D,T)=1$.

For any $D\in \mathcal{C}_{2k+1}^1$ for some $k\ge 2$, by Claim~\ref{claim:C_{2k+1}^1}, we let $S_D\subseteq V(D)$ be a set of $\left \lfloor \tfrac{4k+2}{3}  \right \rfloor $ vertices such that $V(D)-S_D$ is a set of isolated vertices.

For any $D\in \mathcal{C}_{2k+1}^2$ for some $k\ge 2$, we let $S_D\subseteq V(D)$ be the set of $2k +1$ vertices such that $e_G(D-S_D,T)=0$. Note that $D-S_D$ is a graph with at least one vertex.

We let $T^{\prime}$ be the set of all these vertices from $T$ that are adjacent in $G$ to a representative vertex of components from $\mathcal{C}_3$. 
As there are $2m$ representative vertices, we have $|T^{\prime}|=m$.

Let $S^{\prime }=S\cup \left (\bigcup _{D\in \bigcup_{k\ge 1}\mathcal{C}_{2k+1}}S_D \right )\cup T^{\prime }$. Then we have 
\begin{eqnarray*}
    |S^{\prime }| 
        & = & |S|+2c_3+\sum_{k\ge 2}\left \lfloor \frac{4k+2}{3} \right \rfloor c_{2k+1}^1+\sum_{k\ge 2}(2k+1)c_{2k+1}^2+m,\\ 
    c(G-S^{\prime})
        & \ge & |T|-|T^{\prime}|+2|T^{\prime}|+\sum_{k\ge 2}c_{2k+1}^2\\
        & = & |S|+\sum_{k\ge 2}kc_{2k+1}^1+c_3+1+m+\sum_{k\ge 2}(k+1)c_{2k+1}^2,
\end{eqnarray*}
where $|T|=|S|+c_3+\sum_{k\ge 2}kc_{2k+1}^1+\sum_{k\ge 2}kc_{2k+1}^2+1$ by \eqref{eq:T-size}.
By Claim~\ref{claim:S-size}, we have $c(G-S^{\prime })\ge 4$.

Then as $\frac{|S|+a}{|S|+b}$ is decreasing in $|S|$ provided $\frac{a}{b}\ge 1$ and $\frac{5c_3/2+a}{3c_3/2 +b}$ is increasing in $c_3$ provided $\frac{a}{b}\le \frac{5}{3}$ for any positive real numbers $a$ and $b$, by Claim \ref{claim:c3-size} and Claim \ref{claim:S-size}, we have
\begin{eqnarray*}
    \frac{|S^{\prime }|}{c(G-S^{\prime })} 
        & \le & \frac{|S|+2c_3+\sum_{k\ge 2}\left \lfloor \frac{4k+2}{3}  \right \rfloor c_{2k+1}^1+\sum_{k\ge 2}(2k+1)c_{2k+1}^2+m}{|S|+\sum_{k\ge 2}kc_{2k+1}^1+c_3+1+m+\sum_{k\ge 2}(k+1)c_{2k+1}^2} \\
        & \le & \frac{|S|+2c_3+\sum_{k\ge 2}\left \lfloor \frac{4k+2}{3}  \right \rfloor c_{2k+1}^1+\sum_{k\ge 2}(2k+1)c_{2k+1}^2+m}{|S|+\sum_{k\ge 2}kc_{2k+1}^1+c_3+m+\sum_{k\ge 2}(k+1)c_{2k+1}^2} \\
        & \le & \frac{\frac{5}{2}c_3+\sum_{k\ge 2}\left (\frac{8}{3}k-\frac{1}{6} \right ) c_{2k+1}^1+\sum_{k\ge 2}\left (\frac{10}{3}k+\frac{1}{6}\right ) c_{2k+1}^2+m+\frac{|T_1|}{6} +\frac{5}{6}p+\frac{q}{2}+3}
        {\frac{3}{2}c_3+\sum_{k\ge 2}\left (\frac{7}{3}k-\frac{5}{6}\right ) c_{2k+1}^1+\sum_{k\ge 2} \left(\frac{7}{3}k+\frac{1}{6}\right ) c_{2k+1}^2+m+\frac{|T_1|}{6} +\frac{5}{6}p+\frac{q}{2}+3}\\ 
        & \le & \frac{6m+\frac{10}{3}p+3q+|T_1|+\frac{14}{3}+\sum_{k\ge 2}\left(\frac{13}{3}k+\frac{2}{3}\right ) c_{2k+1}^1+\sum_{k\ge 2}\left (5k+1\right) c_{2k+1}^2}{4m+\frac{7}{3}p+2q+\frac{2}{3}|T_1|+4+\sum_{k\ge 2}\left (\frac{10}{3}k-\frac{1}{3} \right) c_{2k+1}^1+\sum_{k\ge 2}\left (\frac{10}{3}k+\frac{2}{3}\right) c_{2k+1}^2} \\
        & < & \frac{3}{2},
\end{eqnarray*}
a contradiction to $\tau(G)\ge \frac{3}{2}$. 
The proof of Theorem~\ref{thm} is now complete. 
\qed

\section*{Acknowledgment}
Songling Shan was partially supported by NSF grant DMS-2345869. 
This work was supported by NSFC (No. 12371356).


\begin{thebibliography}{99}

    \bibitem{Bauer2000} D. Bauer, H. J. Broersma, H. J. Veldman, Not every 2-tough graph is Hamiltonian, Discrete Appl. Math., 99 (1--3): 317--321, 2000.

    \bibitem{Bauer2006} D. Bauer, H. J. Broersma, E. Schmeichel, Toughness in graphs-a survey, Graphs Combin., 22: 1--35, 2006.
        
    \bibitem{Berge1958} C. Berge, Sur le couplage maximum d'un graphe, C. R. Acad. Sci. Paris, 247: 258--259, 1958.    
        
    \bibitem{Chvatal1973} V. Chv\'{a}tal, Tough graphs and Hamiltonian circuits, Discrete Math., 5 (3): 215--228, 1973.        
        
  \bibitem{Kanno2019} J. Kanno, S. Shan, Vizing's 2-factor conjecture involving toughness and maximum degree conditions, Electron. J. Combin., 26 (2): P2.17, 2019.
    
    \bibitem{Owens1999} P. J. Owens, Non-Hamiltonian maximal planar graphs with high toughness,  Tatra Mt. Math. Publ., 18: 89--103, 1999.
                
    \bibitem{Shan2024} S. Shan, A construction of a $\frac{3}{2}$-tough plane triangulation with no 2-factor, J. Graph Theory, 109 (1): 5--18, 2025. 
             
     \bibitem{Tutte1947} W. T. Tutte, The factorization of linear graphs, J. Lond. Math. Soc., s1--22 (2): 107--111, 1947.
           
    \bibitem{Tutte1952} W. T. Tutte, The factors of graphs, Canad. J. Math., 4: 314--328, 1952.
              
    \bibitem{Tutte1956} W. T. Tutte, A theorem on planar graphs, Trans. Amer. Math. Soc., 82 (1): 99--116, 1956.
    
    \bibitem{West2011}  D. B. West, A short proof of the Berge-Tutte formula and the Gallai-Edmonds structure theorem, European J. Combin., 32 (5): 674--676, 2011.


\end{thebibliography}
\end{document}